\documentclass[11pt,a4paper]{article}

\usepackage[left=3.5cm, right=3.5cm, top=1in, bottom=1in]{geometry}

\usepackage{authblk}
\usepackage{graphicx}
\usepackage{color}
\usepackage{amssymb}
\usepackage{psfrag}
\usepackage{xcolor}

\usepackage{subcaption}
\usepackage{fancyhdr}

\begin{document}

% --- TITRE PRINCIPAL ---
\title{Isotropy and Galilean invariance of Lattice Boltzmann Method: Theoretical and numerical analysis using oblique dipole benchmark\thanks{Contribution submitted for publication in {\textit{Computers \& Fluids}} and presented at the ICMMES Conference, 2025.}}

% --- BLOC DES AUTEURS ET AFFILIATIONS ---
\author[1,2]{Fran\c{c}ois Dubois\thanks{fdubois@lecnam.fr}}
\author[3]{Souleyman Kadri-Harouna\thanks{souleyman.kadri-harouna@univ-lr.fr}}
\author[4]{Pierre Lallemand\thanks{plalleman1@free.fr}}
\author[3]{Mohamed Mahdi Tekitek\thanks{Corresponding author: mohamed.tekitek@univ-lr.fr}}

\affil[1]{\small Conservatoire National des Arts et M{\'e}tiers, Laboratoire de M{\'e}canique des Structures et des Syst{\`e}mes Coupl{\'e}s, F-75003 Paris, France}
\affil[2]{\small Laboratoire de Math{\'e}matiques d’Orsay, CNRS, Universit{\'e} Paris-Saclay, B{\^a}t. 307, F-91405 Orsay, France}
\affil[3]{\small MIA Laboratory, La Rochelle Universit{\'e}, Ave. Albert Einstein, BP 33 060, 17 031 La Rochelle, France}
\affil[4]{\small Beijing Computational Science Research Center, Zhangguancun Software Park II, Haidian District, Beijing 100094, China}

% --- GÉNÉRATION DU TITRE ---
\maketitle

% --- FORCE LA PREMIÈRE PAGE À RESTER VIERGE D'EN-TÊTE ---
\thispagestyle{plain}

% --- RÉSUMÉ ET MOTS-CLÉS ---
\begin{abstract}
This work focuses on the two-dimensional, nine-velocity (D2Q9) lattice Boltzmann model.
First, we show that the D2Q9 scheme cannot achieve second-order accuracy unless the cubic
velocity terms are neglected, and we explain how some of these parasitic terms can be
eliminated. Second, we demonstrate that the standard choice of the equilibrium
distribution has no effect on the equivalent PDE at second order. Finally, we numerically
investigate the effect of these cubic terms and study different choices of equilibrium
distributions using a new benchmark called the Oblique Dipole Benchmark, which describes
obliquely propagating 2D vortex dipoles with periodic boundary conditions.

\vspace{0.5cm}
\noindent \textbf{Keywords:} LB scheme with projection, Cubic Terms, MRT, Numerical dissipation, Taylor expansion, ABCD method, Isotropy, Galilean invariance, double vortex
\end{abstract}

% =========================================================================
% DÉBUT DE VOTRE ARTICLE ICI (ex: \section{Introduction})
% =========================================================================

%%%%%%%%%% \linenumbers

%%%%%%%%%%%%%%%%%%%%%%%%%%%%%%%%%%%%%%%%%%%%%% section 0
\section*{Introduction}
Dubois and Lallemand \cite{DL22} show that the choice of the equilibrium distribution function is central to recovering the Navier-Stokes equations at second order in the lattice Boltzmann framework. However, residual second-order defects persist in D2Q9 and D3Q19 schemes, primarily due to cubic terms. This issue was highlighted by Dellar \cite{Del14} and others \cite{QZL98, HJ06a, HJ06b}. It raises questions regarding equilibrium construction in BGK, multiple-relaxation-time (MRT), and elaborate schemes like the projected D2Q9 \cite{DL22}. These residual defects impact Galilean invariance, which must be assessed when the effective shear viscosity is independent of the advection velocity.

The choice of equilibrium deserves particular attention. This choice has no effect on the second-order equivalent PDE, so it can be selected freely at that level. However, it does impact the third order, as previously emphasized by Qian {\it{et al.}} \cite{QHL92} and recently by Lallemand et al. \cite{LDL24}. This impact is crucial because it can fundamentally affect stability and physical modeling. To assess these effects, we designed a numerical test case. This test reveals dispersive errors associated with odd higher-order terms. It also shows the effective numerical viscosity induced by even higher-order truncation errors, particularly regarding their impact on isotropy and Galilean invariance.

The proposed benchmark, inspired by Bruneau and Clercx \cite{BC06}, consists of a double-vortex (dipole) flow with periodic boundary conditions instead of no-slip (Dirichlet) conditions. Unlike standard tests such as the Taylor-Green vortex or lid-driven cavity, this case involves vortical structures advected by a spatially non-uniform velocity field. Consequently, errors stemming from cubic terms and higher-order numerical viscosity (often negligible in classical configurations) become pronounced enough for quantitative analysis. In this benchmark, the double-vortex is self-advected by a non-constant velocity field. The isotropy of the scheme is easily measured by the loss of symmetry or the deviation from a pseudo-spectral reference solution. A further advantage is that the dipole can be advected at any desired angle relative to the mesh with a chosen Reynolds number.
Notably, while \cite{DL22, QZL98, HJ06a, HJ06b} and \cite{LDL24} utilize advected decaying shear or Gaussian stream function flows, those tests focus primarily on Galilean invariance (cubic defects). In such cases, many higher-order terms are canceled because the advection velocity is constant.

In this work, we briefly introduce the Lattice Boltzmann Method (LBM), then present the different choices of equilibrium distributions and study the second-order equivalent PDE for  D2Q9 model.
We focus particularly on MRT and projected models with specific choices of relaxation rates.
Then, we perform our new benchmark and study the influence of different equilibrium choices, considering not only the second-order accuracy but also higher-order corrections.

%%%%%%%%%%%%%%%%%%%%%%%%%%%%%%%%%%%%%%%%%%%%%% section 1
\section{Lattice Boltzmann scheme and equivalent PDE}

\noindent The Lattice Boltzmann Method (LBM) does not start from the partial differential equations (PDEs) to be simulated (unlike classical methods such as finite differences or finite volumes). Instead, it is derived from a discretization of the continuous Boltzmann equation. Specifically, it replaces the continuous velocity space with a finite set of $q$ discrete velocities and discretizes space and time on a regular lattice with parameters $\Delta x$ and $\Delta t$, respectively. Thus a scale speed $\lambda$ is defined by $\lambda=\Delta x/\Delta t$. In a DdQq lattice Boltzmann model,
the resulting numerical scheme describes the evolution of the distribution function $f(x,v_i,t)=f_i(x,t)$ at the time $t+\Delta t$ as follows:
\begin{equation}
f_{i}({x},t+\Delta t)=f_{i}^{*}({x}-v_{i}\Delta t,t), \quad 0\leq i
 \leq q-1,
 \label{lbm_f}
 \end{equation}
where $v_i=\lambda e_i$ are the discrete velocities and the  the superscript $*$ denotes post-collision quantities.
alternates between a transport step, which is exact along lattice characteristics, and a collision step, which models local relaxation toward an equilibrium state.
Among the most widely used collision models are the Bhatnagar–Gross–Krook (BGK) formulation\cite{QHL92}, based on a single relaxation time, and the Multiple-Relaxation-Time (MRT)\cite{DDH92} formulation, in which different moments relax at distinct rates.
We consider here the MRT formulation called "Generalized lattice Boltzmann equations"\cite{DDH92}.
Let $M$ be a constant invertible matrix. Using this matrix, the moment vector is defined as a linear function of the $f_i$ by
\begin{equation}
m_k = \sum_j M_{kj} f_j .
\label{mom}
\end{equation}
The matrix $M$ is constructed by selecting homogeneous polynomials to define the moments,
followed by the application of the Gram--Schmidt algorithm to obtain an orthogonal matrix $M$ (see \cite{DL22} for details).
Let $N$ denote the number of macroscopic partial differential equations (PDE) to be modeled.
The first $N$ moments of the vector $W = (m_1,\ldots,m_N)$ are conserved during the collision step.
The remaining moments, which are not conserved, are gathered in the vector
$
Y = (m_{N+1},\ldots,m_{q-N}).
$
The relaxation step is then described by
\begin{equation}
m^*=\left( \begin{array}{c} W^* \\ Y^* \end{array} \right) = \left(\begin{array}{c} W \\ Y + S (\Phi(W) -Y) \end{array} \right)
\label{col}
\end{equation}
where $\Phi(W) = (m^{eq}_{N+1},\ldots,m^{eq}_{q-N})= Y^{eq}$ denotes the equilibrium distribution,
and $S = \mathrm{diag}(s_1,s_2,\ldots,s_{q-N})$ is a diagonal relaxation matrix of order $q-N$.
Let us recall that the equilibrium distribution $f^{eq}$ can be recovered as
$$f^{eq} = M^{-1} \left(\begin{array}{c} W \\ Y^{eq} \end{array} \right) = M^{-1} \left(\begin{array}{c} W \\ \Phi(W) \end{array} \right)$$
Since the equilibrium distribution depends only on the conserved variables, the mapping $\Phi$ is uniquely determined by these variables and can be constructed accordingly.
Note that if $s_1 = s_2 = \cdots = s_{q-N} = \frac{1}{\tau}$, the classical BGK model with a single relaxation time is recovered.
Now, in order to understand which PDE is solved by the MRT scheme described by equation \ref{lbm_f}, with a collision operator in moment space as given in equation \ref{col}, we aim to derive the equivalent partial differential equations (also called the modified equations). In the lattice Boltzmann community, this is classically achieved using the Chapman--Enskog expansion. However, in this work we use the ABCD method, introduced by one of the author \cite{Dub22}, which is based on a Taylor expansion.
This approach has the advantage of providing the equivalent PDE without introducing multiscale expansions, auxiliary variables, or additional operators.
To this end, we first rewrite the scheme in moment space. We obtain:
\begin{eqnarray*}
m_k(x,t+\Delta t) &=& \sum_j M_{kj} f_j^{\star}(x - v_j \Delta t, t) \\
&= &\sum_{\ell} \sum_{n=0}^{\infty} \frac{1}{n!}
\left(
\sum_j M_{kj}
\left(-\Delta t \sum_{\alpha} v_j^{\alpha} \partial_{\alpha}\right)^n
(M^{-1})_{j\ell}
\right)
m_\ell^{\star}(x,t) \\
&= &\bigl(\exp(-\Delta t\,\Lambda)\, m^{\star}(x,t)\bigr)_k.
\end{eqnarray*}
where the advection operator $\Lambda$ is defined by
$
\Lambda = M \, \mathrm{diag}\!\left(\sum_{\alpha=1}^{d} v^{\alpha}\partial_{\alpha}\right) M^{-1}.
$
%The exponential operator can be expanded as
%$
%\exp(-\Delta t\,\Lambda)
%= I - \Delta t\,\Lambda + \frac{\Delta t^2}{2}\Lambda^2
%+ \cdots + (-1)^n \frac{\Delta t^n}{n!}\Lambda^n + \cdots .
%$
We consider the following block decomposition of the advection operator:
$$
\Lambda \equiv \left( \begin{array}{cc}
A & B \\
C & D
\end{array} \right)
$$
Using the splitting of the moment vector into conserved $W$ and non-conserved $Y$, as given in \ref{col}, and using the ABCD method (see \cite{Dub22}), the following equivalent PDE at order $2$ is obtained:
\begin{equation}
\partial_t W + \Gamma_1 + \Delta t\,\Gamma_2 = O(\Delta t^2),
\label{EDP2}
\end{equation}
where
\begin{eqnarray*}
\Gamma_1 &=& A\,W + B\,\Phi(W), \\
\Gamma_2 &=& B\,\Sigma\,\Psi_1,\\
\Psi_1 &=& D\Phi(W)\,\Gamma_1 - \bigl(C\,W + D\,\Phi(W)\bigr).
\end{eqnarray*}
Where $\Sigma$ is the Hénon matrix given by $\Sigma \equiv S^{-1} - \frac{1}{2}$.\\
Note that the second-order equivalent PDE \ref{EDP2} is obtained under the assumption of a constant velocity scale $\lambda =\frac{\Delta x}{\Delta t}$, corresponding to the acoustic scaling. The present second-order expansions are derived in full generality, {\it{i.e}}., for an arbitrary spatial dimension $d$ and an arbitrary discrete velocity set of cardinality $q$.\\
Remark: It is shown in \cite{Dub22} that the ABCD methodology yields the same macroscopic PDE as those obtained via the classical Chapman--Enskog asymptotic expansion, up to fourth order accuracy.\\
Without loss of generality, we henceforth restrict the analysis to the
two-dimen-\\sional case and consider nine discrete velocities, denoted as D2Q9.
The number of conserved moments is $N = 3$, corresponding to the density $\rho$ and the momentum components $\rho u$ and $\rho v$ in the $x$ and $y$ directions, respectively, in order to model the isothermal Navier-Stokes equations. If only a single conserved moment $\rho$ is retained ($N = 1$), one recovers the advection--diffusion equation \cite{MTB22,MTB23}.
\subsection*{The D2Q9 scheme}
\begin{figure}[htbp!]
\begin{center}
\includegraphics[width=4cm,angle=0]{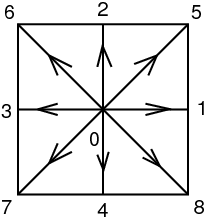}
\end{center}
\caption{Stencil for the D2Q9 lattice Boltzmann scheme.}
\label{stensild2q9}
\end{figure}

\noindent We consider the classical D2Q9 lattice Boltzmann scheme\cite{LL00}. Let $\mathcal{L}^0$ denote the regular lattice of nodes
$\mathcal{L}^0 \equiv \left\{ x_{i,j} \in (\Delta x \mathbb{Z}, \Delta y \mathbb{Z}) \right\},$
with uniform spacing $\Delta x = \Delta y$. The model is based on nine discrete velocities
$v_i = \lambda\, e_i, \ i = 0,\ldots,8,$
corresponding to the D2Q9 velocity set. The discrete velocity directions are defined as
$e_0 = (0,0),$ $e_1 = (1,0),$ $e_2 = (0,1),$ $e_3 = (-1,0),$ $e_4 = (0,-1),$ $e_5 = (1,1),$ $e_6 = (-1,1),$ $e_7 = (-1,-1),$ $e_8 = (1,-1),$
as illustrated in Figure~\ref{stensild2q9}.
%%%
In this case the d'Humieres matrix of moments is given by:
\begin{equation} {{M}}=\left ( \begin{array}{ccccccccc}
\quad 1 &\quad 1 &\quad 1&\quad 1&\quad 1&\quad 1&\quad 1 &\quad 1&\quad 1\\
\quad 0
&\lambda & 0&-\lambda&0&\lambda&-\lambda&-\lambda&\lambda \\ \quad 0 &0&
\lambda & 0&-\lambda&\lambda&\lambda&-\lambda&-\lambda\\
  -4\lambda^2&-\lambda^2&-\lambda^2&-\lambda^2&-\lambda^2&2\lambda^2&2\lambda^2&2\lambda^2&2\lambda^2\\
  0 &\lambda^2 &-\lambda^2&\lambda^2&-\lambda^2&0&0&0&0\\
  0 &0 & 0&0&0&\lambda^2&-\lambda^2&\lambda^2&-\lambda^2 \\
  0 &-2 \lambda^3&0&2\lambda^3&0&\lambda^3&-\lambda^3&-\lambda^3&\lambda^3\\
  0 & 0&-2\lambda^3 & 0&2\lambda^3&\lambda^3&-\lambda^3&-\lambda^3&\lambda^3\\
  4\lambda^4 &-2\lambda^4&-2\lambda^4&-2\lambda^4&-2\lambda^4&\lambda^4&\lambda^4&\lambda^4&\lambda^4
 \end{array} \right).
\label{matrixq9}
\end{equation}
Using the linear transformation given by equation \ref{mom}, we get the following moments:
$m_0= \rho$ density, $m_1=j_x=\rho u$, $m_2=j_y=\rho v$ are momentum components in $x$ and $y$ directions.
$m_3=\varepsilon$ is energy, $m_4=xx$, $m_5=xy$  are diagonal and off-diagonal component of the stress tensor.
$m_6=q_x$, $m_7=q_y$ are  the x and y components of energy flux and $m_8=h$ is the square energy.

\noindent {\bf{Remark :}} The d’Humières moment matrix $M$ is orthogonal. The different moments are obtained by considering polynomial functions of the velocity components $(u,v)$: zeroth-order polynomials correspond to the density, first-order polynomials to the momentum components $j_x$ and $j_y$, second-order polynomials to the energy $\varepsilon$ and the stress components $xx$ and $xy$, third-order polynomials to the heat flux components $q_x$ and $q_y$, and fourth-order polynomials to the higher-order moment $h$.

\noindent In order to apply the ABCD method, we identify the conserved moments
$W = (\rho, \rho u, \rho v),$
and the non-conserved moments $Y = (\varepsilon, xx, xy, q_x, q_y, h).$
The corresponding equilibrium distribution functions are then defined as follows:
\begin{equation}
Y^{eq} = \Phi(W) =
\left( \begin{array}{c}
\Phi_\varepsilon \\
\Phi_{xx}\\
\Phi_{xy}\\
\Phi_{qx}\\
\Phi_{qy}\\
\Phi_{h}
\end{array}
\right) =
\left( \begin{array}{c}
m_3^{eq}\\
m_4^{eq}\\
m_5^{eq}\\
m_6^{eq}\\
m_7^{eq}\\
m_8^{eq}
\end{array} \right) =
 \left( \begin{array}{c}
-2 \rho +3 \rho {\lambda^2} (u^2+v^{2}) \\
\frac{1}{\lambda^2}\left(u^2-v^2\right) \\
 \rho \frac{1}{\lambda^2} u v \\
 -\frac{\rho u}{\lambda} \\
 -\frac{\rho v}{\lambda} \\
 \rho - 3 \rho {\lambda^2} (u^2+v^{2})
\end{array} \right)
\label{equi_qian}
\end{equation}
Note that the above equilibrium moments correspond to the following equilibrium distribution function $f^{{eq}}$:
$$f_{j}^{eq}=t_{j}\rho\left[1+\frac{3}{\lambda}(e_{j}.{\bf{u}})+\frac{9}{2\lambda^{2}}(e_{j}.{\bf{u}})^2-\frac{3}{2\lambda^{2}}|{\bf{u}}|^{2}\right], \qquad j=0..8,$$

\noindent where ${\bf{u}} = (u,v)$, $t_0 = \frac{4}{9}$, $t_{1,2,3,4} = \frac{1}{9}$, and $t_{5,6,7,8} = \frac{1}{36}$.
This is the standard equilibrium distribution used in the lattice Boltzmann method for the isothermal Navier--Stokes equations, originally introduced in the pioneering work~\cite{QHL92}.

\noindent Thus, for the scheme defined by Eqs.~$(\ref{lbm_f})$ and~$(\ref{col})$, with the equilibrium distribution given by~$(\ref{equi_qian})$, the second-order equivalent partial differential equation~$(\ref{EDP2})$, obtained using the ABCD method, reads as follows:

\begin{equation}
\begin{array}{ccc} \displaystyle \frac{\partial \rho}{\partial t}+ \frac{\partial
\rho u}{\partial x}+\frac{\partial \rho v}{\partial y}&=&O(\Delta t^2) ,
\end{array}
\label{edp1}
\end{equation}
\begin{equation}
\begin{array}{rl}
\displaystyle\frac{\partial \rho u}{\partial t} &\displaystyle + \frac{\lambda^2 }{3} \frac{\partial \rho}{\partial x} + \frac{\partial \rho u^2}{\partial x} + \frac{\partial \rho u v}{\partial y} = \nonumber \\
&= \displaystyle  \frac{\lambda^2\Delta t  }{3}\left[ \sigma_e \rho \frac{\partial}{\partial x} \left(\frac{\partial u}{\partial x}+\frac{\partial v}{\partial y}\right) + \sigma_x \rho \Delta u \right]  + \mathcal{S}_x  + O(\Delta t^2)\nonumber
\end{array}
\label{edp2}
\end{equation}
\begin{equation}
\begin{array}{cl}
\displaystyle \frac{\partial \rho v}{\partial t} &\displaystyle + \frac{\lambda^2 }{3} \frac{\partial \rho}{\partial y} + \frac{\partial \rho u v}{\partial x} + \frac{\partial \rho v^2}{\partial y} = \nonumber \\
&= \displaystyle  \frac{\lambda^2\Delta t }{3}\left[ \sigma_e \rho \frac{\partial}{\partial y} \left(\frac{\partial u}{\partial x}+\frac{\partial v}{\partial y}\right) + \sigma_x \rho \Delta v \right] + \mathcal{S}_y+ O(\Delta t^2)\nonumber
\end{array}
\label{edp3}
\end{equation}
%where the relaxation rates are given by
%$s_1 = s_e, \quad s_2 = s_3 = s_x, \quad s_4 = s_5 = s_q, \quad s_6 = s_h$ and $\Sigma = \mathrm{diag}(\sigma_e, \sigma_x, \sigma_x, \sigma_q, \sigma_q, \sigma_h).$
%and \sigma_i=\left(\frac{1}{s_i}-\frac{1}{2}\right).
where the relaxation rates are given by $s_1 = s_e$, $s_2 = s_3 = s_x$,
$s_4 = s_5 = s_q$, and $s_6 = s_h$. The associated matrix $\Sigma$ is defined as
$\Sigma = \mathrm{diag}(\sigma_e, \sigma_x, \sigma_x, \sigma_q, \sigma_q, \sigma_h),$
with the components given by
$\sigma_i = \left(\frac{1}{s_i} - \frac{1}{2}\right).$
Details of the matrices $A$, $B$, $C$, $D$, $\Gamma_1$, and $\Gamma_2$ can be found in~\cite{DL22}.
So at second order, the above resulting equivalent PDE recover the isothermal Navier–Stokes equations, with a sound speed $c_s= \lambda/\sqrt{3}$, a bulk viscosity $\zeta = \frac{\lambda^2\Delta t }{3}\sigma_e$ , and a shear viscosity $\nu=  \frac{\lambda^2\Delta t }{3}\sigma_x$.
The two terms $\mathcal{S}_x$ and $\mathcal{S}_y$ are spurious contributions. In order to recover the exact Navier--Stokes equations at second order with full Galilean invariance, these terms must vanish. However, in the present case of the D2Q9 scheme, these terms are non-zero and are given by:
\begin{eqnarray}
\mathcal{S}_x &=&\frac{\Delta t}{2} \partial_x \Big[ (\sigma_e + \sigma_x) \left( \partial_x (\rho u^3) + \partial_y (\rho u^2 v) \right) +  \label{S_X}\\
&+& (\sigma_e - \sigma_x) \left( \partial_x (\rho u v^2) + \partial_y (\rho v^3) \right) + 2\rho uv (\sigma_e \partial_x v + \sigma_x \partial_y u) \Big]+ \nonumber\\
&+&\Delta t \, \sigma_x \, \partial_y \left[ \partial_x (\rho u^2 v) + \partial_y (\rho u v^2) \right]\nonumber \\
\nonumber \\
\mathcal{S}_y &=& \frac{\Delta t}{2} \partial_y \Big[ (\sigma_e + \sigma_x) \left( \partial_y (\rho v^3) + \partial_x (\rho u v^2) \right) + \label{S_Y}\\
& +& (\sigma_e - \sigma_x) \left( \partial_y (\rho u^2 v) + \partial_x (\rho u^3) \right) + 2\rho uv (\sigma_e \partial_y u + \sigma_y \partial_x v) \Big] +\nonumber\\
& + &\Delta t \, \sigma_x \, \partial_x \left[ \partial_y (\rho u v^2) + \partial_x (\rho u^2 v) \right] \nonumber
\end{eqnarray}

Note that, as discussed in~\cite{DL22}, the choice of equilibrium for the second-order moments
$\varepsilon$, $xx$, and $xy$ (i.e., the moments associated with second-order polynomial invariants)
determines the structure of the equivalent partial differential equations at first order.
These moments are therefore referred to as the \emph{Eulerian moments}, $Y_e = (\varepsilon, xx, xy).$
The moments $q_x$ and $q_y$ control the second-order contributions of the equivalent PDEs, in particular the viscous stress tensor.
Thus, we obtain the first scheme, denoted as $LB1$, where the equilibrium
is given by \ref{equi_qian}. This corresponds to the standard LBM scheme
\cite{QHL92,LL00} for the isothermal Navier-Stokes equations.

Note that the choice of the equilibrium for the highest-order moment $h$ does
not affect the second-order equivalent equations; it therefore has no impact
on $(\ref{edp1})$, $(\ref{edp2})$, and $(\ref{edp3})$. Consequently, the remaining
moments are referred to as the \emph{viscous moments}, defined as $Y_v = (q_x, q_y, h)$.

Hence, in the second scheme, denoted as $LB2$, we employ D2Q9 scheme
to investigate the influence of the equilibrium associated with the square of the energy, $h$.
While the classical equilibrium \ref{equi_qian} prescribes
\begin{equation}
h^{eq} = \rho - 3 \rho \lambda^2 (u^2 + v^2),
\end{equation}
we propose a simplified alternative where the nonlinear contribution is canceled:
\begin{equation}
h^{eq} = \rho.
\end{equation}
This modification is introduced to assess its numerical impact, given that it solely
affects the third- and higher-order terms of the equivalent partial differential equations.
%%%%%%%%
%%% LB3 projeté
\noindent Next, we introduce the third scheme, denoted as $LB3$, which is based on the projected D2Q9 LB scheme proposed by one of the authors in \cite{DP22}. The core principle of this approach is to project the viscous moments, $\Phi_v = (q_x, q_y, h)$, as functions of the conserved moments, $W = (\rho, \rho u, \rho v)$, and the Eulerian moments, $\Phi_e = (\varepsilon, xx, xy)$, according to
\begin{equation}
\Phi_v = K W + L \Phi_e,
\end{equation}
where the projection matrices are defined as
\begin{equation}
K = \left(
\begin{array}{ccc}
0 & -\lambda^2 & 0 \\
0 & 0 & -\lambda^2 \\
-\lambda^4 & 0 & 0
\end{array}\right),
\qquad
L =\left(
\begin{array}{ccc}
0 & 0 & 0 \\
0 & 0 & 0 \\
-\lambda^2 & 0 & 0
\end{array}\right).
\end{equation}
Consequently, the collision step for the viscous moments reduces to
\begin{equation}
q_x^* = -\lambda^2 \rho u, \qquad
q_y^* = -\lambda^2 \rho v, \qquad
h^* = -\lambda^4 \rho - \lambda^2 \varepsilon^* .
\end{equation}
It is worth noting that this projected scheme yields the same equivalent partial differential equations, namely \ref{edp1}, \ref{edp2}, and \ref{edp3}, as the standard MRT scheme defined by \ref{lbm_f} and \ref{col} with the equilibrium distribution \ref{equi_qian}. The key advantages of this projected formulation are twofold: it requires the tuning of only a single relaxation rate, $s_e$, and it exhibits significantly improved numerical stability (see \cite{DP22} for further details).
%%%%
%%.   LB4
%%%
Furthermore, the existence of the spurious terms $\mathcal{S}_x$ and $\mathcal{S}_y$
has already been investigated in \cite{Del14, HJ06a, HJ06b, QZL98}. These terms
are responsible for the appearance of a non-Galilean invariant viscous stress tensor.
More recently, \cite{DL22} proposed an alternative choice for the equilibrium moments,
given by
\begin{eqnarray}
\Phi_{q_x} = q_x^{eq} &=& -\frac{\rho u}{\lambda} + 3\frac{\rho u}{\lambda^2}\left(u^2+v^2\right), \\
\Phi_{q_y} = q_y^{eq} &=& -\frac{\rho v}{\lambda} + 3\frac{\rho v}{\lambda^2}\left(u^2+v^2\right),
\label{du_eq}
\end{eqnarray}
while the same equilibrium as in \ref{equi_qian} is retained for the remaining
non-conserved moments. Note that this equilibrium is non-standard, as noted by Philippi~\cite{Phi18}.
With this specific choice, the spurious terms contain fewer contributions.
In fact, applying the ABCD method yields the following expressions for these terms:
\begin{eqnarray}
\mathcal{S}_x &=& \Delta t \, \sigma_x \left( \partial_x \left[ \partial_x (\rho u^3) - \partial_y (\rho v^3) \right] - \partial_y \left[ \partial_x(\rho v^3) + \partial_y (\rho u^3) \right] \right), \label{S_X_FD} \\
\mathcal{S}_y &= &\Delta t \, \sigma_x \left( \partial_y \left[ \partial_y (\rho v^3) - \partial_x (\rho u^3) \right] - \partial_x \left[ \partial_y (\rho u^3) + \partial_x (\rho v^3) \right] \right). \label{S_Y_FD}
\end{eqnarray}
This formulation fully defines the third scheme, hereafter referred to as $LB4$.

\noindent Finally, we introduce the fifth scheme, denoted as $LB5$, which exploits the fact that the relaxation parameters $s_q$ and $s_h$ do not appear in the second-order equivalent partial differential equations \ref{edp1}--\ref{edp3}. Consequently, their values remain free at this order.
These parameters can therefore be used to improve boundary accuracy \cite{GA94, DLT10}, or to enhance the accuracy of the method \cite{DL09}, or to stabilize the scheme.
In particular, the study by Augier \emph{et al.} \cite{ADGL14} utilizes this additional freedom to improve the isotropy properties of the model.
Following this approach, we perform numerical investigations of these isotropic variants of the D2Q9 scheme.
The corresponding relaxation parameters for $LB5$ scheme are selected as follows:
\begin{equation}
s_\varepsilon = s_x, \quad s_q \mbox{ such that } \sigma_q \sigma_x = \frac{1}{6} \quad \mbox{ and } \quad s_h = s_x.
\end{equation}

\noindent To summarize this section, the following D2Q9 lattice Boltzmann schemes are considered in the numerical experiments:
\begin{table}[h!]
\centering
\caption{Summary of the D2Q9 lattice Boltzmann schemes investigated in this work.}
\label{tab:d2q9_schemes}
\begin{tabular}{ll}
\hline
\textbf{Scheme} & \textbf{Description} \\
\hline
LB$1$ & Standard MRT-LBM with the equilibrium distribution~$(\ref{equi_qian})$ \\
LB$2$ & LBM with modified equilibrium $h^{eq} = \rho$ \\
LB$3$ & Projected LBM scheme \\
LB$4$ & LBM with Dubois equilibrium for $q_x$ and $q_y$ \\
LB$5$ & LBM with Dubois equilibrium with Augier isotropy-enhanced \\
& relaxation parameters \\

\hline
\end{tabular}
\end{table}

\newpage

\section{The oblique dipole benchmark}
\begin{figure}[htbp!]
\begin{center}
\includegraphics[width=0.6\textwidth,angle=0]{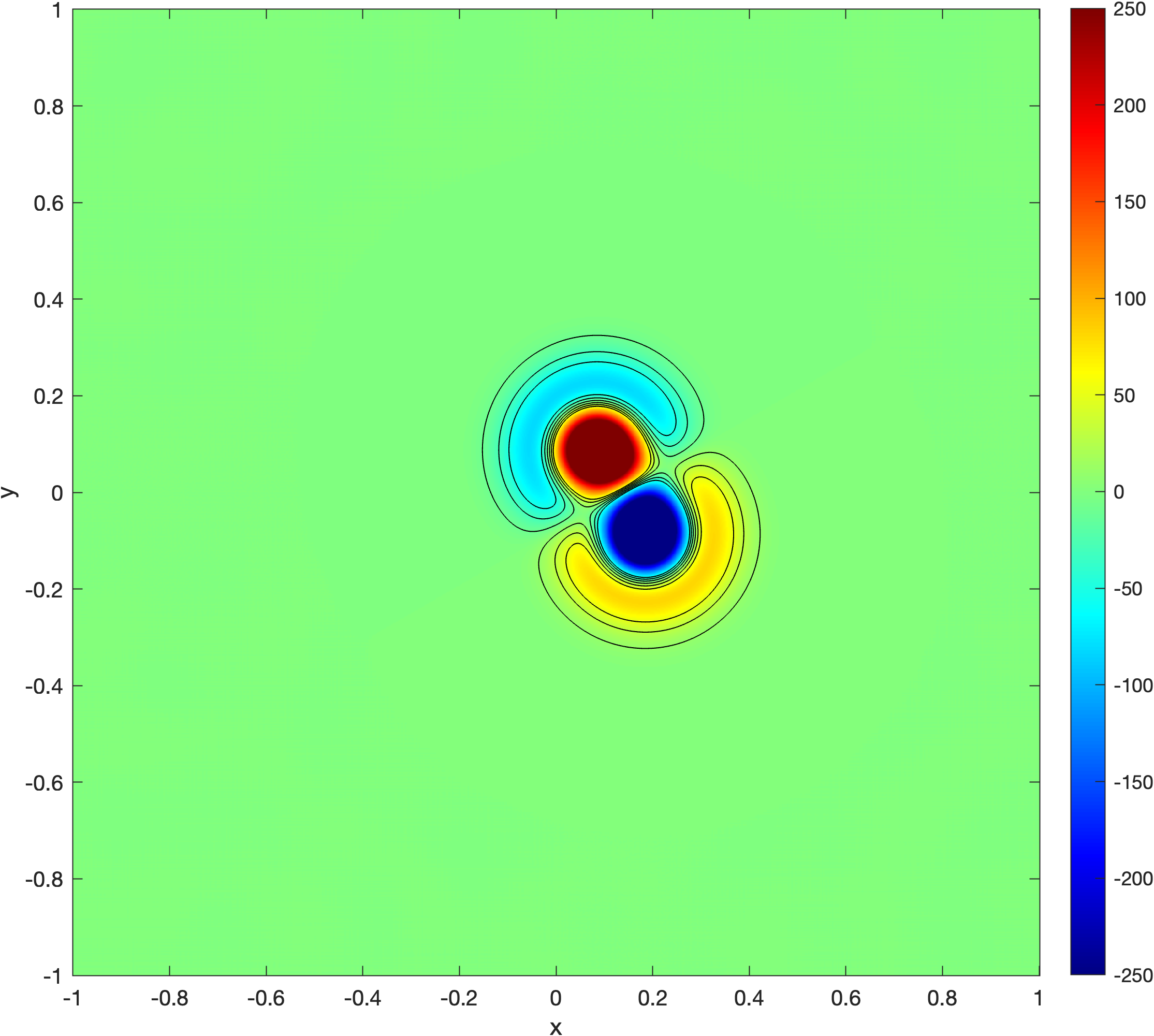}
\end{center}
\caption{Initial vorticity field.}
\label{fig:initial_dipole}
\end{figure}

\noindent A numerical benchmark is introduced to evaluate the isotropy properties of the proposed numerical scheme.
The benchmark is derived from the classical two-dimensional vortex dipole problem, which has been extensively investigated as a prototype flow for studying vortex dynamics and vortex-wall interactions. Early numerical studies were reported by Orlandi \cite{OR90}, while subsequent investigations were conducted by Bruneau {\it et al.} \cite{BC06}. The configuration has also been considered in the context of the lattice Boltzmann method in \cite{LC07,MGR18}.
In the present work, the original benchmark is modified in order to specifically characterize numerical isotropy. While the initial vortex dipole is retained, the solid-wall boundary condition is replaced by periodic boundary conditions in both spatial directions. This modification eliminates the influence of boundary treatments and enables the assessment of anisotropic errors originating exclusively from the spatial discretization.
The objective of this benchmark is not to assess the physical accuracy of the flow solution itself, but rather to quantify the anisotropic errors introduced by the numerical discretization.
By removing the influence of boundary conditions, the proposed configuration isolates the impact of the lattice symmetry and spatial operators on the directional accuracy of the scheme.
The resulting configuration provides a controlled framework for quantifying directional biases introduced by the numerical scheme and for comparing its isotropic properties against the reference solution.
Owing to its spectral accuracy, this method provides a reliable reference against which the isotropic behavior of the numerical scheme can be assessed.
The flow consists of a pair of counter-rotating vortices (vortex dipole) evolving in a square computational domain $[-1,1] \times [-1,1]$ subject to periodic boundary conditions in both spatial directions.
The initial velocity field is defined as
       % \begin{equation*}
    $$\left \{ \begin{array}{lcl}
        u_0(x, y) &= &\frac{\omega_e}{2} \left[ (y - y_1) exp(-\frac{r_2^2}{r_0^2}) - (y - y_2) exp(-\frac{r_1^2}{r_0^2}) \right] \\
        v_0(x, y) &= &\frac{\omega_e}{2} \left[ (x - x_1) exp(-\frac{r_1^2}{r_0^2}) - (x - x_2) exp(-\frac{r_2^2}{r_0^2}) \right]
        \end{array} \right.$$
        %\end{equation*}
where $r_i(x, y) = \sqrt{(x - x_i)^2 + (y - y_i)^2}, \quad i \in \{1, 2\}$, with $r_0 = 0.1$ and $\omega_e = 299.52838$.
The Reynolds number for this benchmark is set to $Re = 2500$, based on the initial dipole velocity and core size, which defines the kinematic viscosity $\nu$ used in the simulations.
The centers of the two Gaussian vortices are located at $(x_1,y_1) = (0.08397460, 0.08660254)$ and $(x_2,y_2) = (0.18397460, -0.08660254)$ to have an angle of propagation $\theta = \frac{\pi}{6}$.
Figure~\ref{fig:initial_dipole} shows the initial vorticity field of the vortex dipole.
%%%%%
%%%
%%%%%%%%%%%%
\noindent With the above initial conditions, the two vortices form a dipole that undergoes a self-induced translational motion.
The direction of propagation can be prescribed through the initial centers of the two Gaussian vortices, $(x_1,y_1)$ and $(x_2,y_2)$, allowing the dipole trajectory to form an arbitrary angle $\theta$ with respect to the mesh.
Unlike classical benchmark problems such as the lid-driven cavity flow or the Taylor--Green vortex, where the dominant flow structures remain stationary with respect to the computational domain,
the present configuration involves the continuous advection of coherent vortical structures across the mesh.
As a consequence, the solution is particularly sensitive to directional errors induced by the spatial discretization,
making this benchmark especially suitable for assessing the isotropy properties of the numerical scheme.
Since no analytical solution is available for the considered configuration, a pseudo-spectral solver is used to generate a reference solution.

\subsection*{Reference solution}

\noindent The reference solution is computed using a Fourier pseudo-spectral method for the spatial discretization \cite{Canuto1988,Canuto2006,Peyret2002,Rogallo1981}. After projection onto the frequency domain, the corresponding generic equation is as follows:
\begin{equation}%\label{eqfourier1}
\partial_t \widehat{\mathbf{u}}
=
-\nu |\mathbf{k}|^2 \widehat{\mathbf{u}}
+
\mathbf{H}(\widehat{\mathbf{u}}),
\label{eq:spectral_ns}
\end{equation}

where
\begin{equation}
\mathbf{H}(\widehat{\mathbf{u}})=\widehat{\mathbb{P}[f-(\mathbf{u}\cdot\nabla)\mathbf{u}]}
\end{equation}
with $\mathbb{P}$ the Leray projector
$$
 \widehat{\mathbb{P}(\mathbf{u})}=\widehat{\mathbf{u}} - \frac{\mathbf{k}(\mathbf{k}\cdot\widehat{\mathbf{u}})}{|\mathbf{k}|^2}.
$$
This allows to decouple de computation of the velocity from the pressure. The pressure is recovered through Helmholtz decomposition:
$$
\nabla p=f-(\mathbf{u}\cdot\nabla)\mathbf{u}-\mathbb{P}[f-(\mathbf{u}\cdot\nabla)\mathbf{\mathbf{u}}].
$$
\noindent Equation $(\ref{eq:spectral_ns})$  is an ordinary differential equation for each frequency, and an integrating-factor formulation is used to handle the viscous term precisely. Meanwhile, the nonlinear contribution is handled using a third-order Adams–Bashforth scheme \cite{durran91}. To avoid aliasing or filtering techniques, we used a fourth-order finite difference method at grid points for the spatial derivatives of the nonlinear term. The resulting time integration reads as:

\begin{equation}
\widehat{\mathbf{u}}^{\,n+1}
=
e^{-\nu |\mathbf{k}|^2 \Delta t}
\left(
\widehat{\mathbf{u}}^{\,n}
+
\frac{\Delta t}{12}
\left[
23\,\mathbf{H}(\widehat{\mathbf{u}}^{n})
-
16\,\mathbf{H}(\widehat{\mathbf{u}}^{n-1})
+
5\,\mathbf{H}(\widehat{\mathbf{u}}^{n-2})
\right]
\right).
\label{eq:ab3_if}
\end{equation}

\begin{figure}[htbp!]
\begin{center}
\includegraphics[width=0.8\textwidth,angle=0]{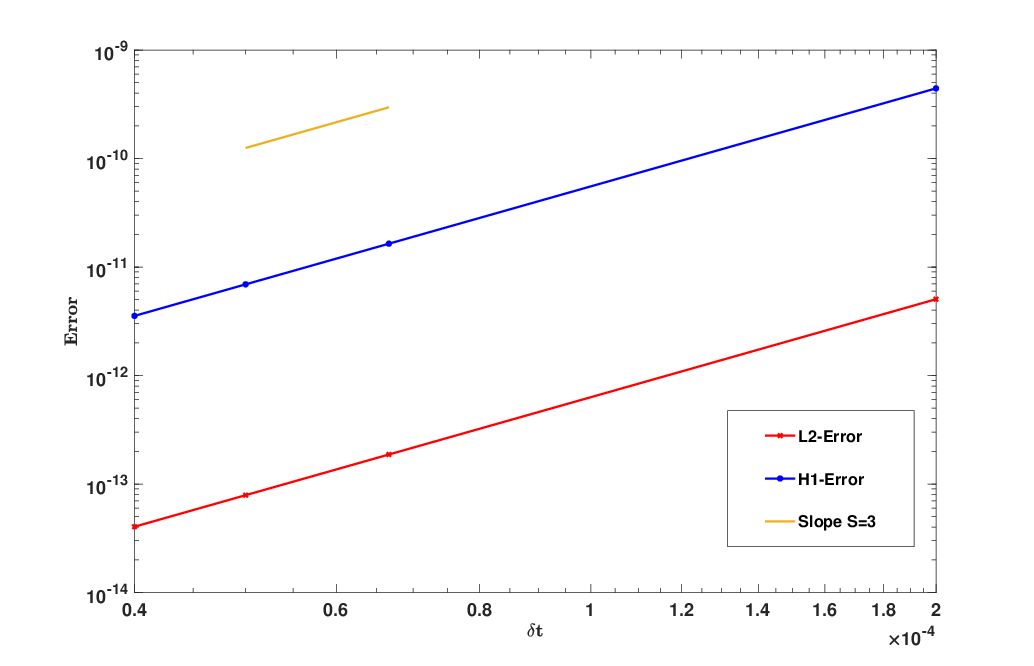}
\end{center}
\caption{Convergence rate of the pseudo-spectral method for the two-dimensional Taylor--Green vortex.}
\label{cv_tg_spec}
\end{figure}

\begin{figure}[htbp]
    \centering
    % =================== PREMI√àRE LIGNE (t = 0.2 √† 0.6) ===================
    \begin{minipage}[b]{0.31\textwidth}
        \centering
        \includegraphics[width=\textwidth]{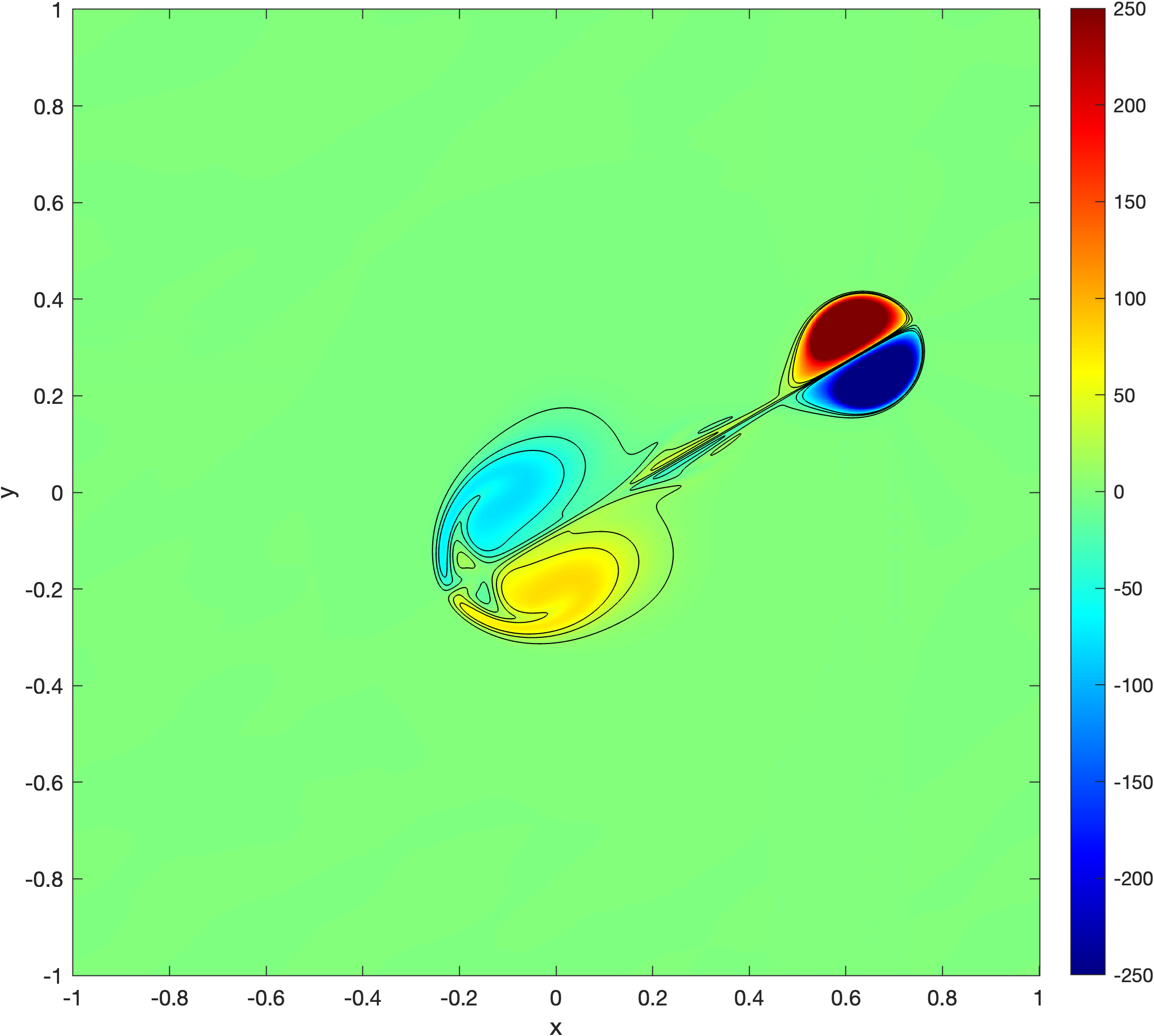}
        \\ \small $t = 0.2$
    \end{minipage}\hfill
    \begin{minipage}[b]{0.31\textwidth}
        \centering
        \includegraphics[width=\textwidth]{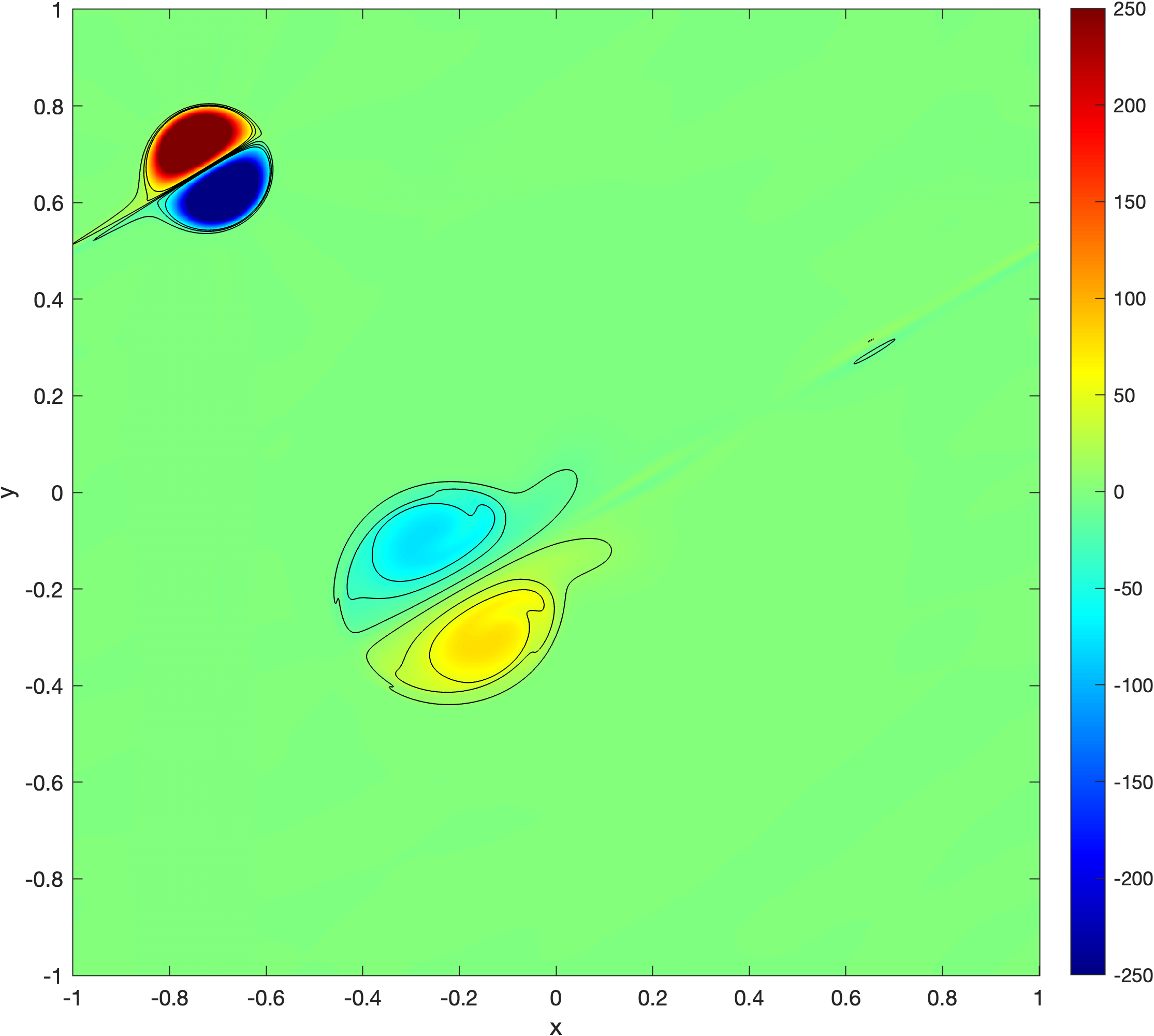}
        \\ \small $t = 0.4$
    \end{minipage}\hfill
    \begin{minipage}[b]{0.31\textwidth}
        \centering
        \includegraphics[width=\textwidth]{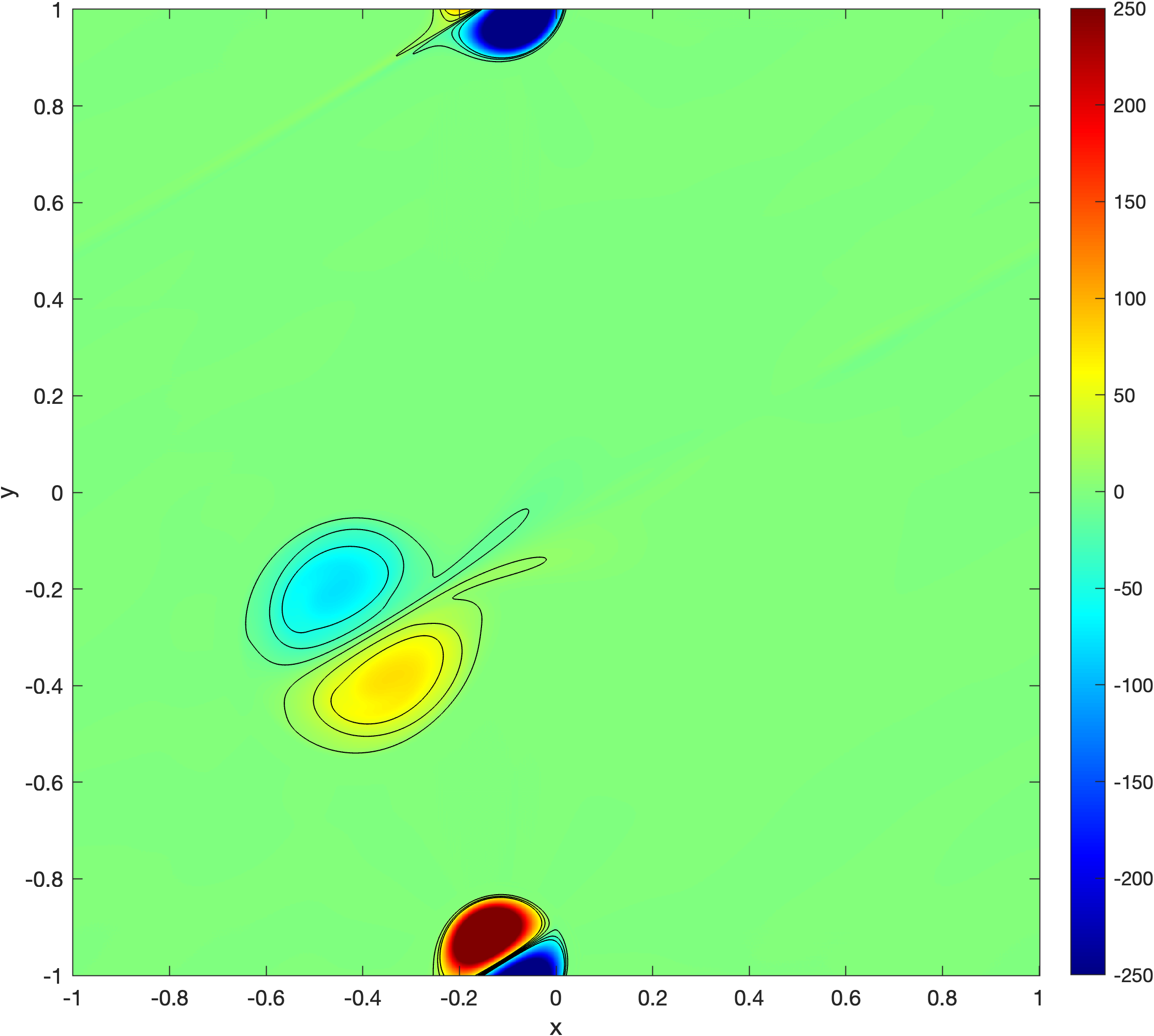}
        \\ \small $t = 0.6$
    \end{minipage}

    \vspace{0.4cm} % Espace vertical entre la ligne 1 et la ligne 2

    % =================== DEUXI√àME LIGNE (t = 0.8 √† 1.2) ===================
    \begin{minipage}[b]{0.31\textwidth}
        \centering
        \includegraphics[width=\textwidth]{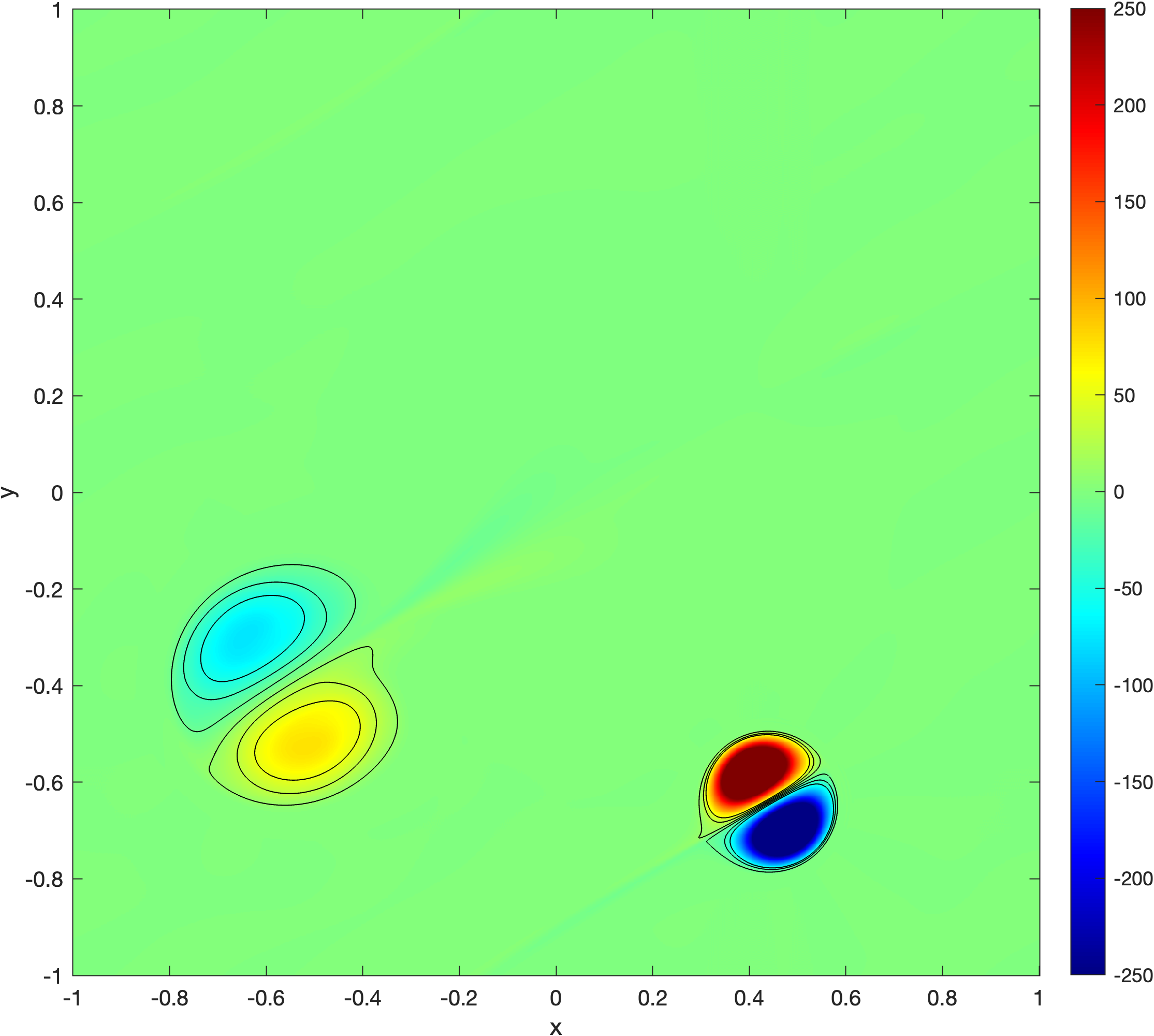}
        \\ \small $t = 0.8$
    \end{minipage}\hfill
    \begin{minipage}[b]{0.31\textwidth}
        \centering
        \includegraphics[width=\textwidth]{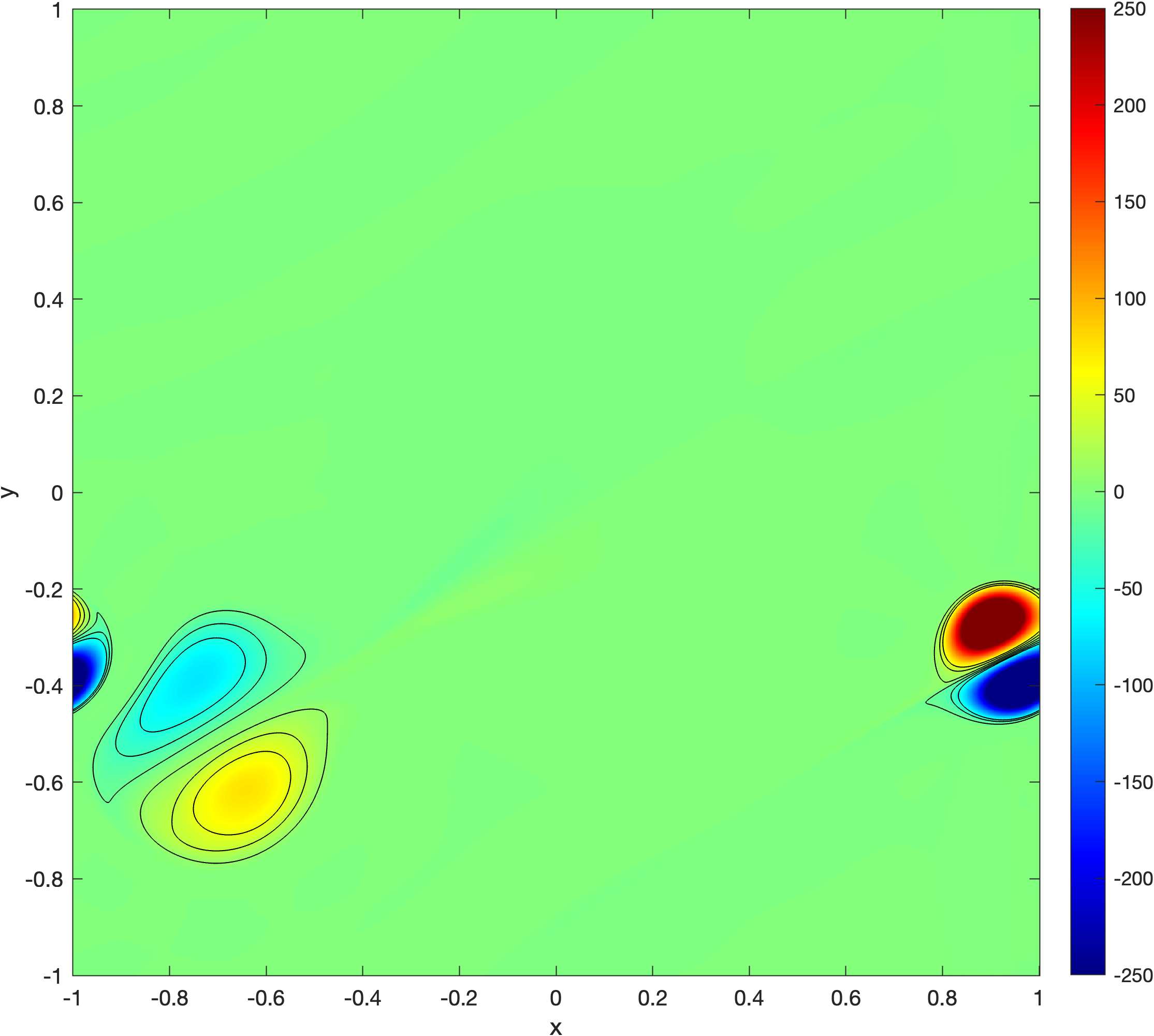}
        \\ \small $t = 1.0$
    \end{minipage}\hfill
    \begin{minipage}[b]{0.31\textwidth}
        \centering
        \includegraphics[width=\textwidth]{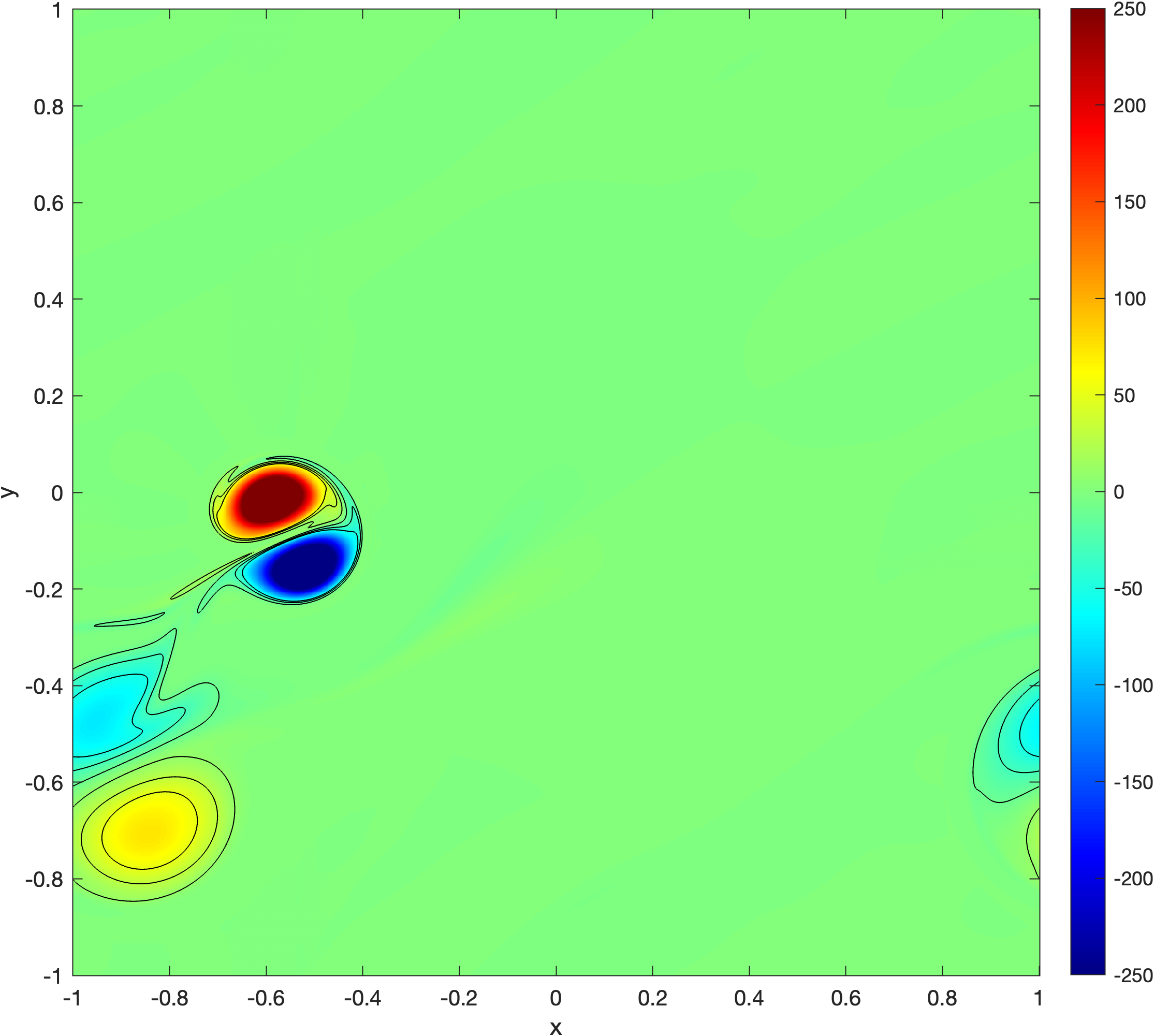}
        \\ \small $t = 1.2$
    \end{minipage}

    \vspace{0.4cm} % Espace vertical entre la ligne 2 et la ligne 3

    % =================== TROISI√àME LIGNE (t = 1.4 √† 1.8) ===================
    \begin{minipage}[b]{0.31\textwidth}
        \centering
        \includegraphics[width=\textwidth]{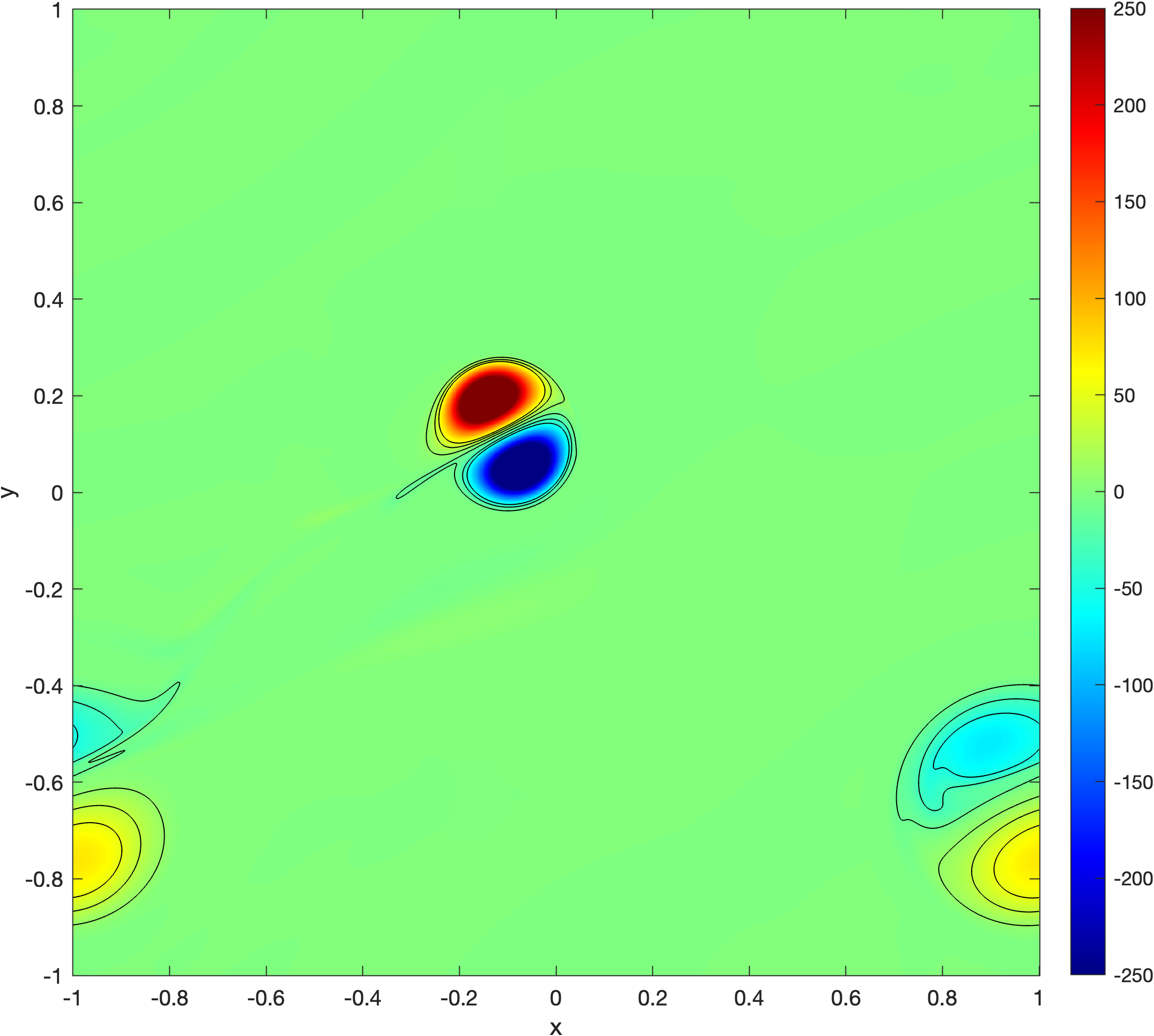}
        \\ \small $t = 1.4$
    \end{minipage}\hfill
    \begin{minipage}[b]{0.31\textwidth}
        \centering
        \includegraphics[width=\textwidth]{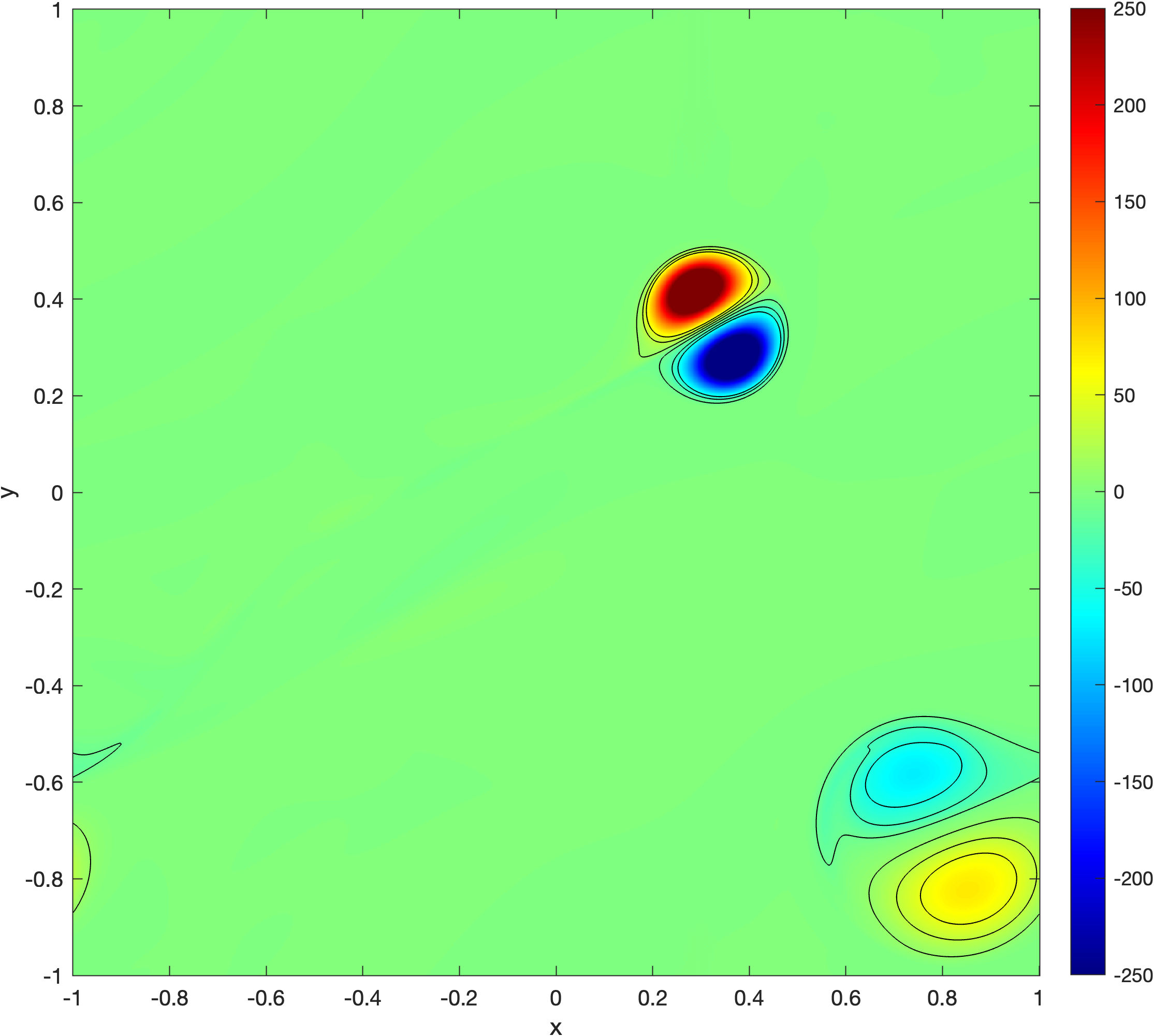}
        \\ \small $t = 1.6$
    \end{minipage}\hfill
    \begin{minipage}[b]{0.31\textwidth}
        \centering
        \includegraphics[width=\textwidth]{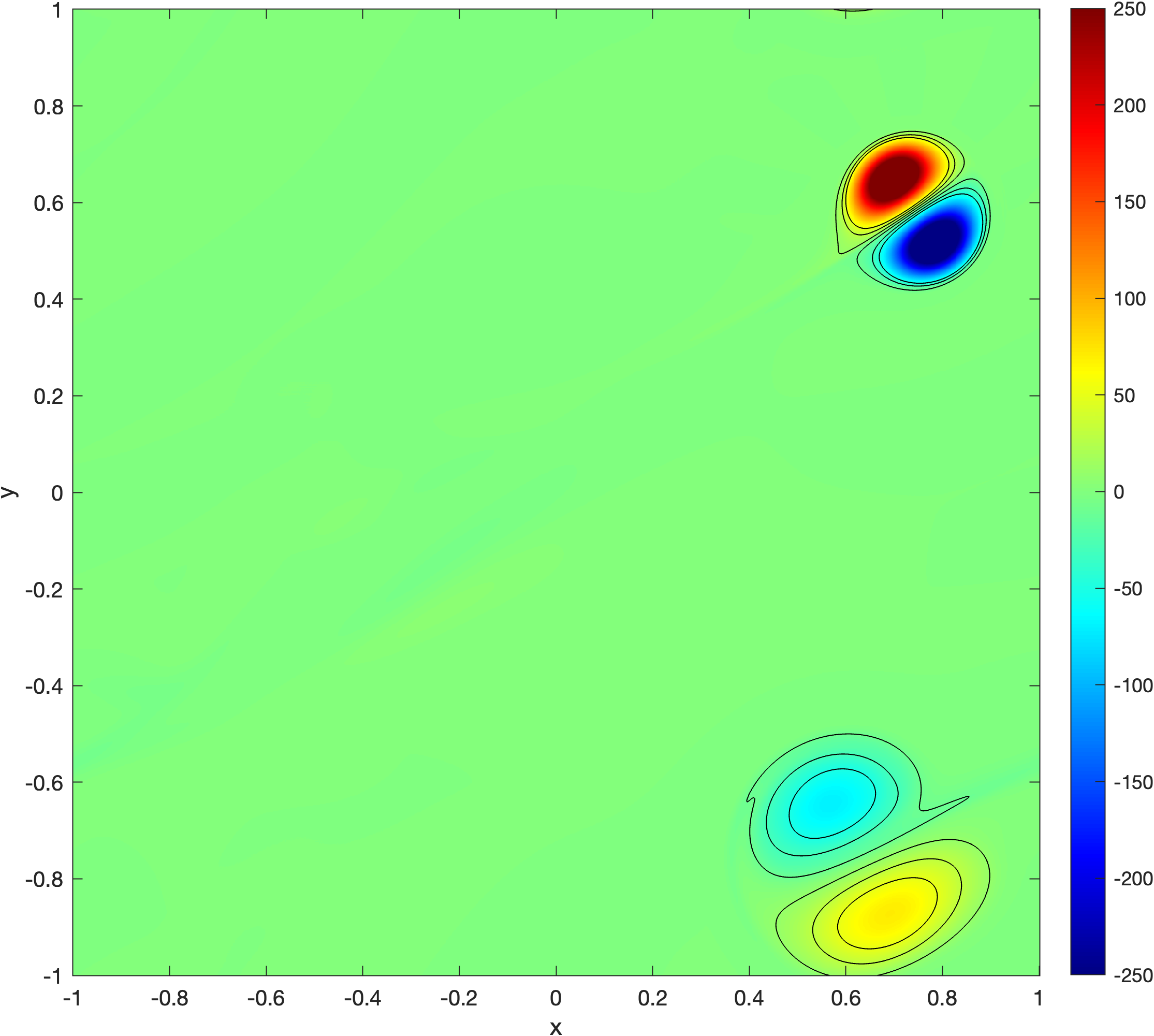}
        \\ \small $t = 1.8$
    \end{minipage}

    % L√©gende globale en anglais pour votre article
    \caption{Reference solution: Time evolution of the vorticity field obtained with the pseudo-spectral solver.}
    \label{fig:vorticite_champ_temporel}
\end{figure}

\noindent Finally, this formulation achieves spectral accuracy in space for diffusion, while also providing third-order accuracy in time and avoiding the stability constraints associated with the diffusive operator. The accuracy of the implementation was verified using the
two-dimensional Taylor--Green vortex problem confirming the expected third-order
convergence rate (see Figure \ref{cv_tg_spec}).
%%%%%%%%%%%%
%%%%%
\noindent Unless otherwise stated, the reference solution is computed on a uniform grid of
$1024^2$ points using a time step $\Delta t = 10^{-4}$. At this resolution,
both spatial and temporal discretization errors are negligible compared with
those of the lattice Boltzmann simulations considered in the present study.
The resulting pseudo-spectral solution is therefore used as the reference
solution for the isotropy assessment.\\
To illustrate the flow evolution obtained with the pseudo-spectral solver, a sequence of vorticity contour plots is presented in Fig.~\ref{fig:vorticite_champ_temporel} at different times,
$t=0.2,\,0.4,\,\ldots,\,2.0$.

\begin{figure}[htbp]
    \centering
    % --- Premi√®re figure : √ânergie ---
    \begin{minipage}[b]{0.32\textwidth}
        \centering
        \includegraphics[width=\textwidth]{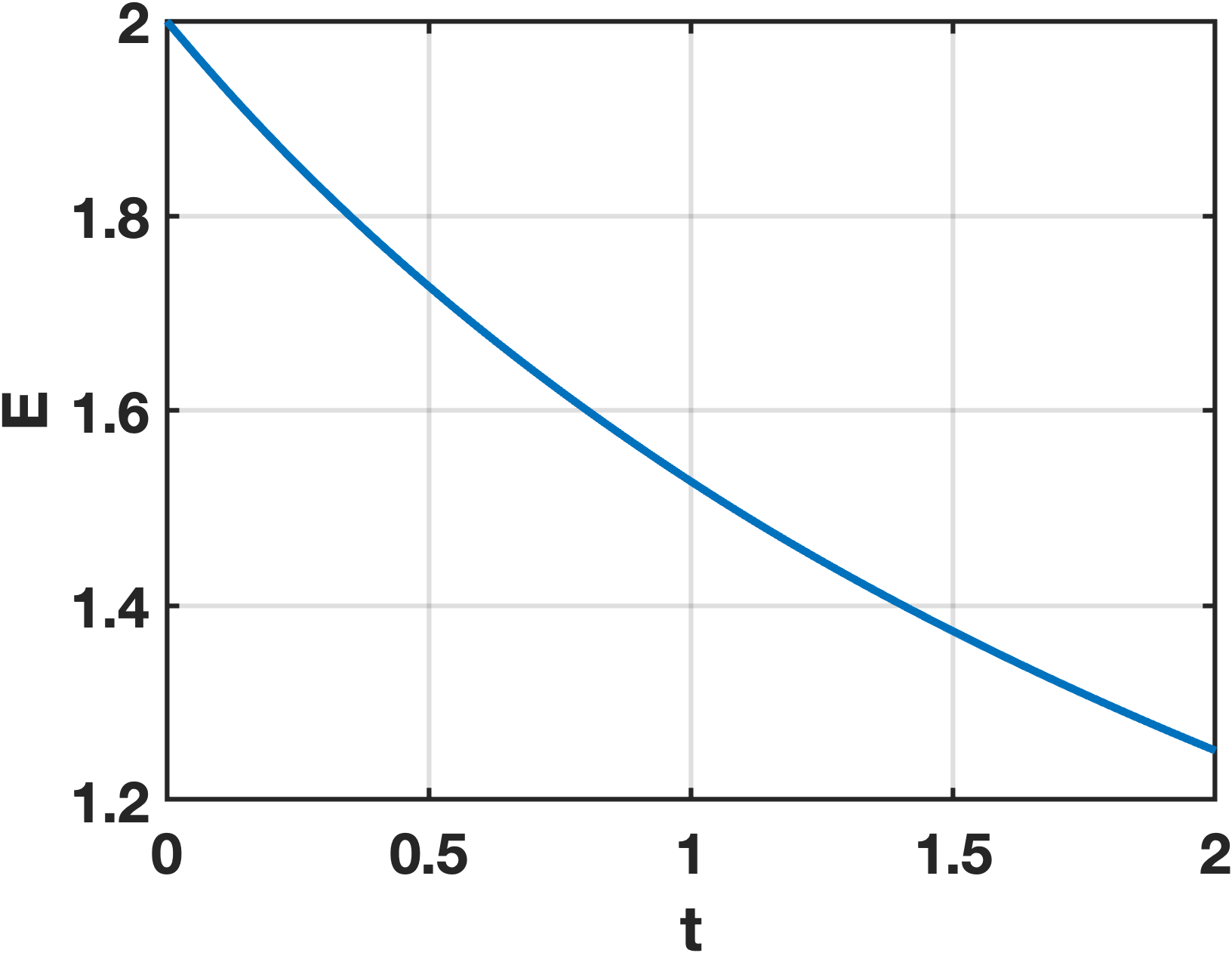}
        \\ \small (a) $E(t)$
    \end{minipage}
    \hfill % Espace horizontal entre les images
    % --- Deuxi√®me figure : Enstrophie ---
    \begin{minipage}[b]{0.32\textwidth}
        \centering
        \includegraphics[width=\textwidth]{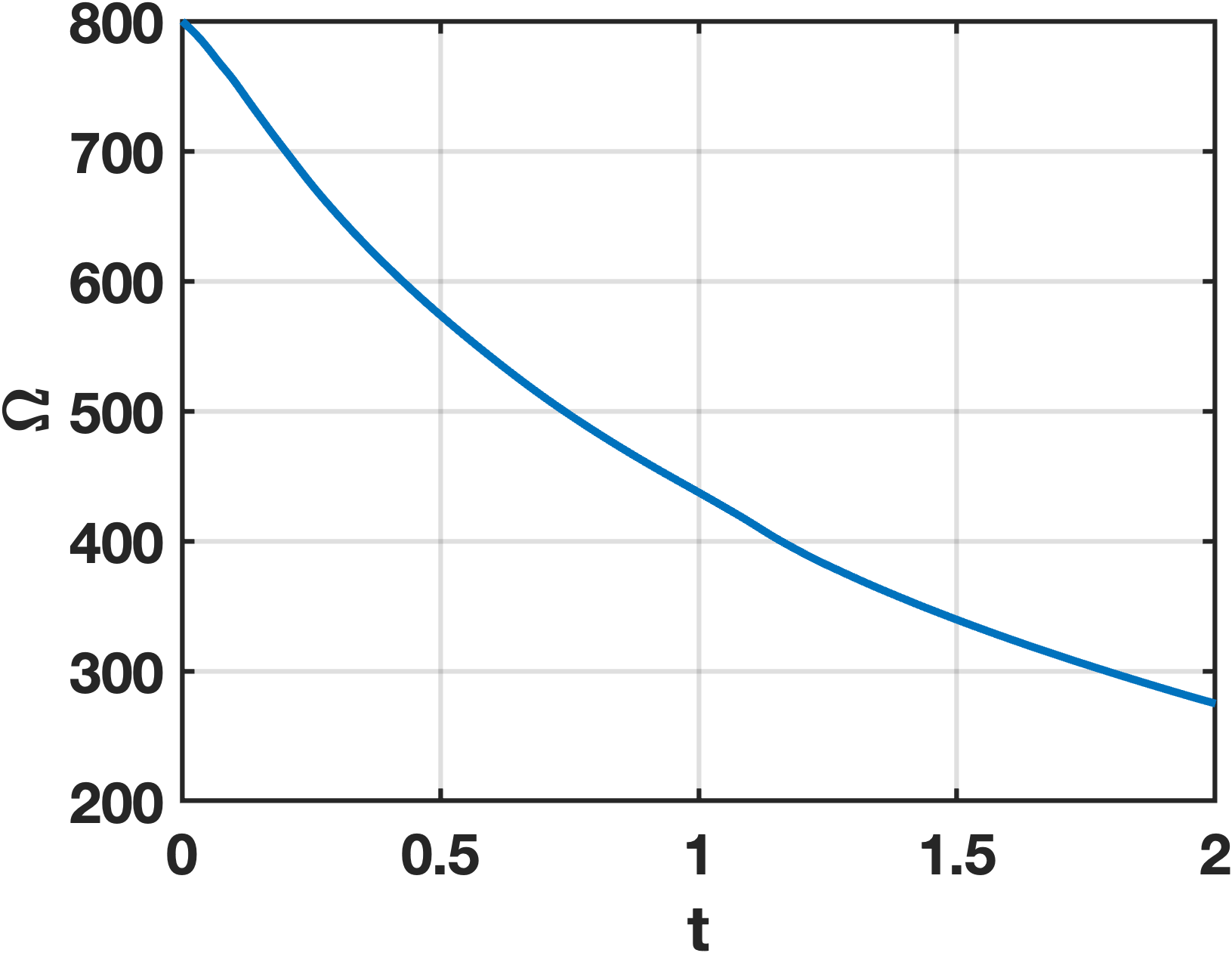}
        \\ \small (b) $\Omega(t)$
    \end{minipage}
    \hfill % Espace horizontal entre les images
    % --- Troisi√®me figure : Palinstrophie ---
    \begin{minipage}[b]{0.32\textwidth}
        \centering
        \includegraphics[width=\textwidth]{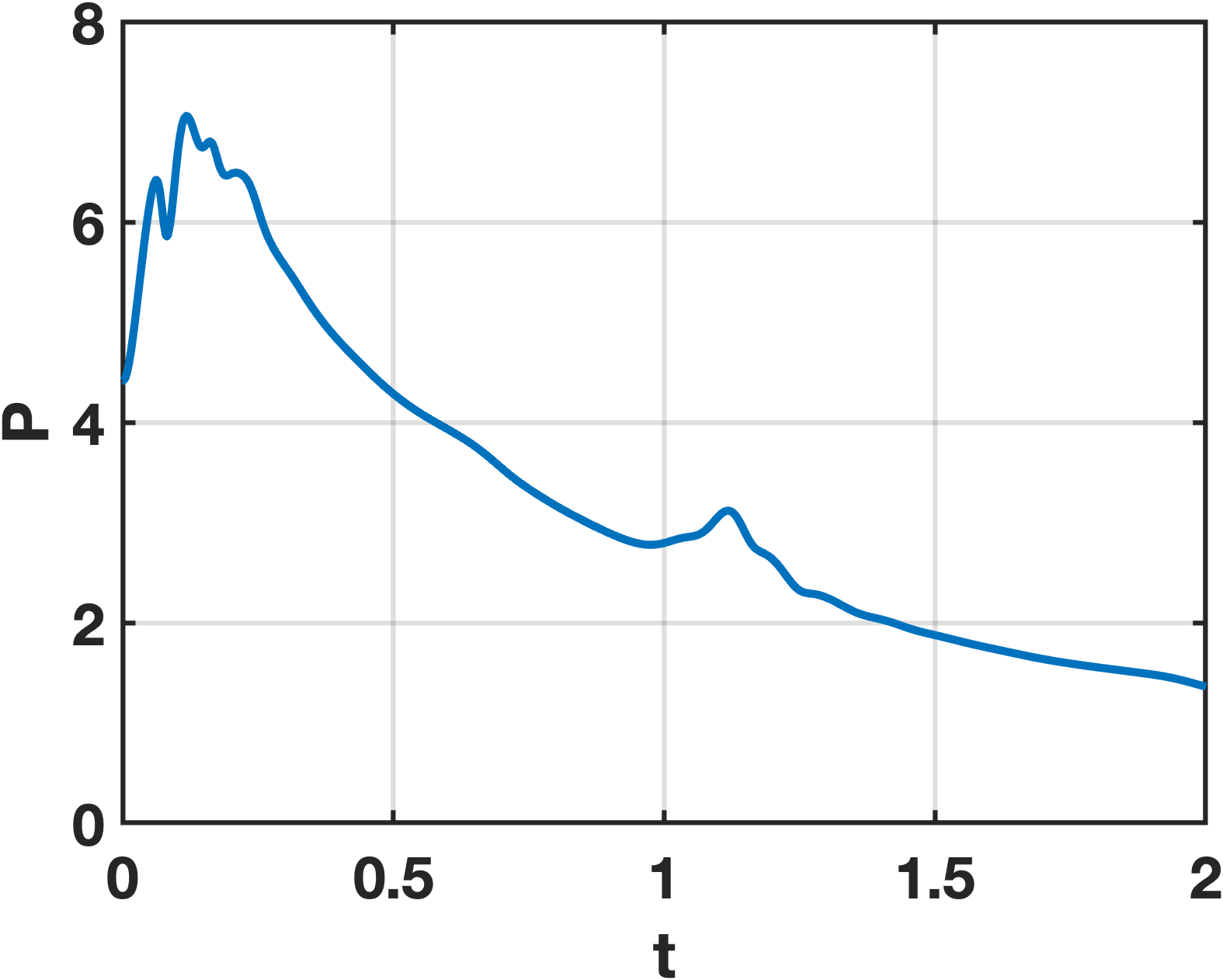}
        \\ \small (c) $P(t)$
    \end{minipage}

    % L√©gende g√©n√©rale unique pour l'ensemble de la ligne
    \caption{Reference solution: Time evolution of the energy (a), enstrophy (b), and palinstrophy (c) computed using the pseudo-spectral method.}
    \label{fig:global_quantities}
\end{figure}

\noindent For a quantitative comparison with the lattice Boltzmann simulations, three global quantities are monitored throughout the computation: the total kinetic energy \(E(t)\), the enstrophy \(\Omega(t)\), and the palinstrophy \(P(t)\). The latter provides a measure of the intensity of vorticity gradients within the flow and is therefore particularly sensitive to small-scale structures and numerical dissipation.
The total kinetic energy is defined as
\begin{equation}
E(t)
=
\frac{1}{2}
\int_{-1}^{1}
\int_{-1}^{1}
\|\mathbf{u}(\mathbf{x},t)\|^2
\,dx\,dy,
\label{eq:energy}
\end{equation}
while the enstrophy is given by
\begin{equation}
\Omega(t)
=
\frac{1}{2}
\int_{-1}^{1}
\int_{-1}^{1}
\omega^2(\mathbf{x},t)
\,dx\,dy,
\label{eq:enstrophy}
\end{equation}
where
\begin{equation}
\omega
=
\frac{\partial u_y}{\partial x}
-
\frac{\partial u_x}{\partial y}
\label{eq:vorticity}
\end{equation}
denotes the scalar vorticity.
Finally, the palinstrophy is defined as
\begin{equation}
P(t)
=
\frac{1}{2}
\int_{-1}^{1}
\int_{-1}^{1}
\left|
\nabla \omega(\mathbf{x},t)
\right|^2
\,dx\,dy.
\label{eq:palinstrophy}
\end{equation}
These integral quantities provide complementary information on the flow dynamics and constitute relevant indicators for assessing the accuracy and isotropy properties of the numerical schemes under consideration.
Note that the vortex strength parameter \(\omega_e=  299.52838\) is chosen such that the initial kinetic energy satisfies $E(0)=2.$
Figure~\ref{fig:global_quantities} shows the time evolution of the kinetic energy \(E\), the enstrophy \(\Omega\), and the palinstrophy \(P\) obtained with the pseudo-spectral solver.

\section{Numerical solution with different LB schemes models}
\noindent In this section, we numerically investigate the different LB schemes presented in Table~\ref{tab:d2q9_schemes} by simulating the oblique dipole benchmark described above, and compare the results with the reference solution.
For all simulations, a uniform grid resolution of $2048 \times 2048$ is employed with $\Delta x = \Delta t$. The spatial domain is defined as $[-1,1] \times [-1,1]$, and the flow is simulated over the time interval $t \in [0, 2]$. Naturally, the Reynolds number is kept identical to that of the reference solution, namely $Re = 2500$. For each scheme, the time evolution of the global kinetic energy $(\ref{eq:energy})$, enstrophy $(\ref{eq:enstrophy})$ and palinstrophy $(\ref{eq:palinstrophy})$ obtained by the LB simulation is compared against the reference solution.
Furthermore, the vorticity field computed by the LB scheme at time $t = 1.8$ is provided, and the positions of the two main vortical structures are compared with those of the reference solution.

%%%%%%%%%
%%%%%%%%%
%%%%%%%%%.      LB1
%%%%%%%%%
%%%%%%%%%

\subsection*{LB1 scheme}
\begin{figure}[htbp]
    \centering
    % --- Première figure : Comparaison Énergie ---
    \begin{minipage}[b]{0.32\textwidth}
        \centering
        \includegraphics[width=\textwidth]{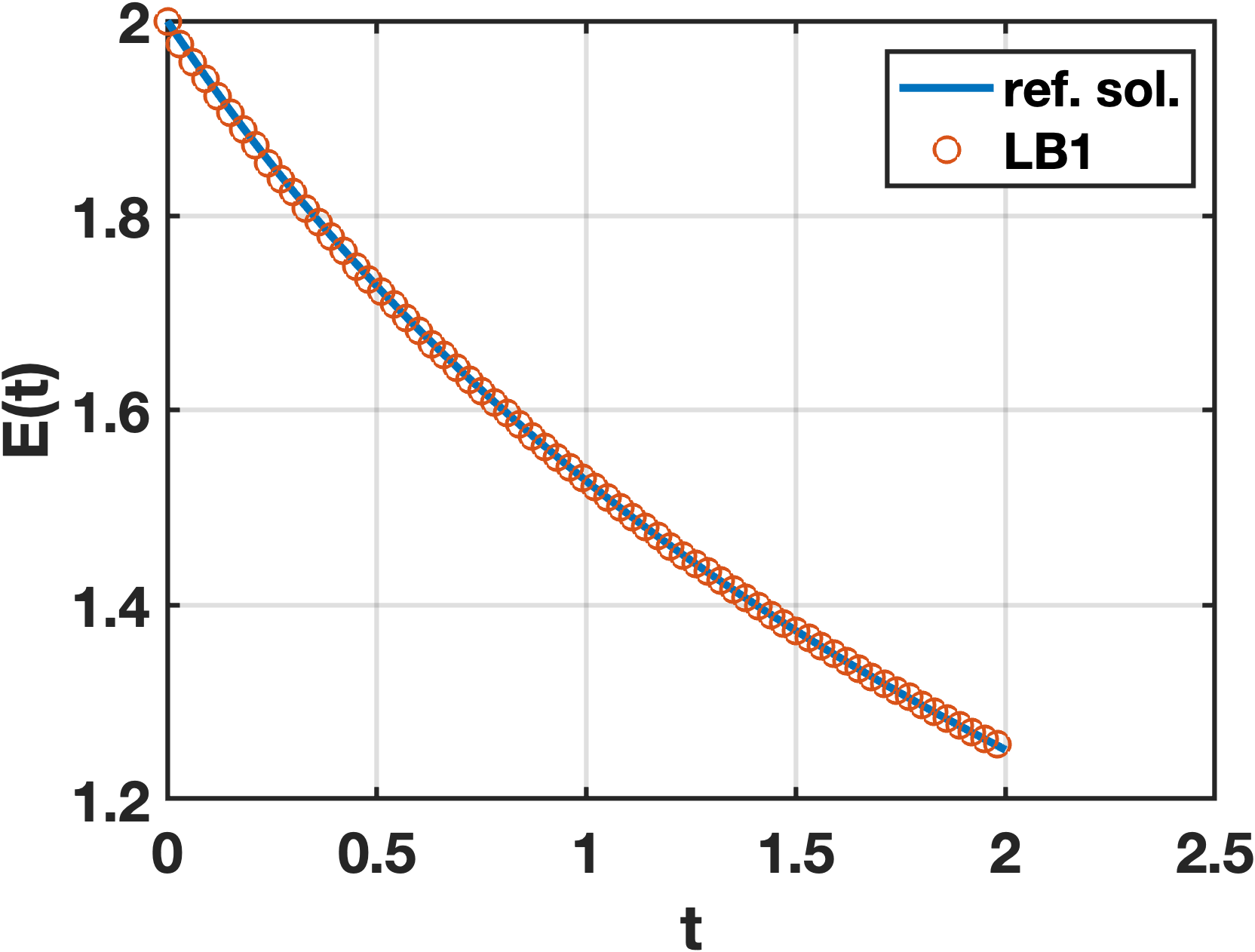}
        \\ \small (a) Kinetic energy $E(t)$
    \end{minipage}
    \hfill % Ajuste l'espace horizontal de manière homogène
    % --- Deuxième figure : Comparaison Enstrophie ---
    \begin{minipage}[b]{0.32\textwidth}
        \centering
        \includegraphics[width=\textwidth]{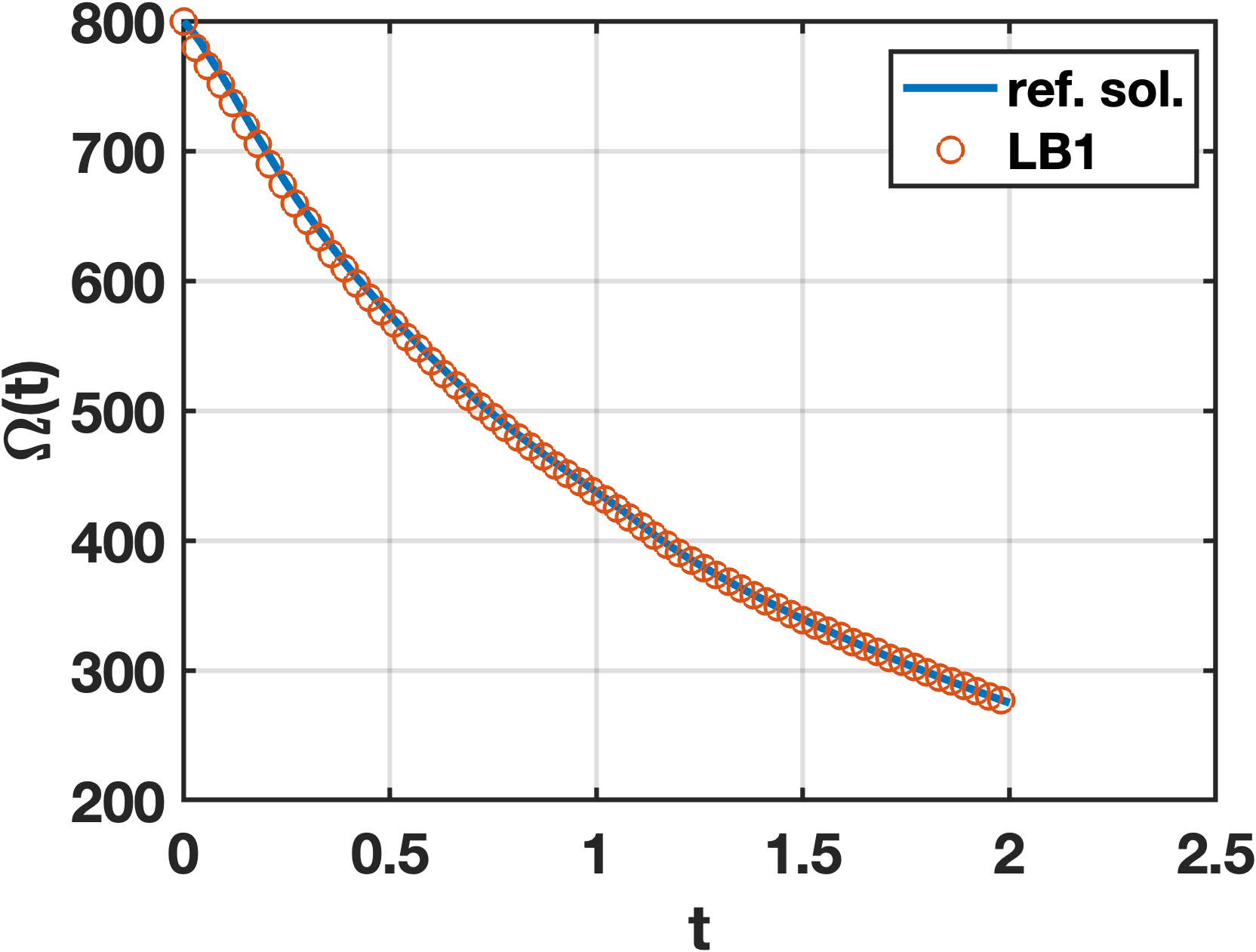}
        \\ \small (b) Enstrophy $\Omega(t)$
    \end{minipage}
    \hfill % Ajuste l'espace horizontal de manière homogène
    % --- Troisième figure : Comparaison Palinstrophie ---
    \begin{minipage}[b]{0.32\textwidth}
        \centering
        \includegraphics[width=\textwidth]{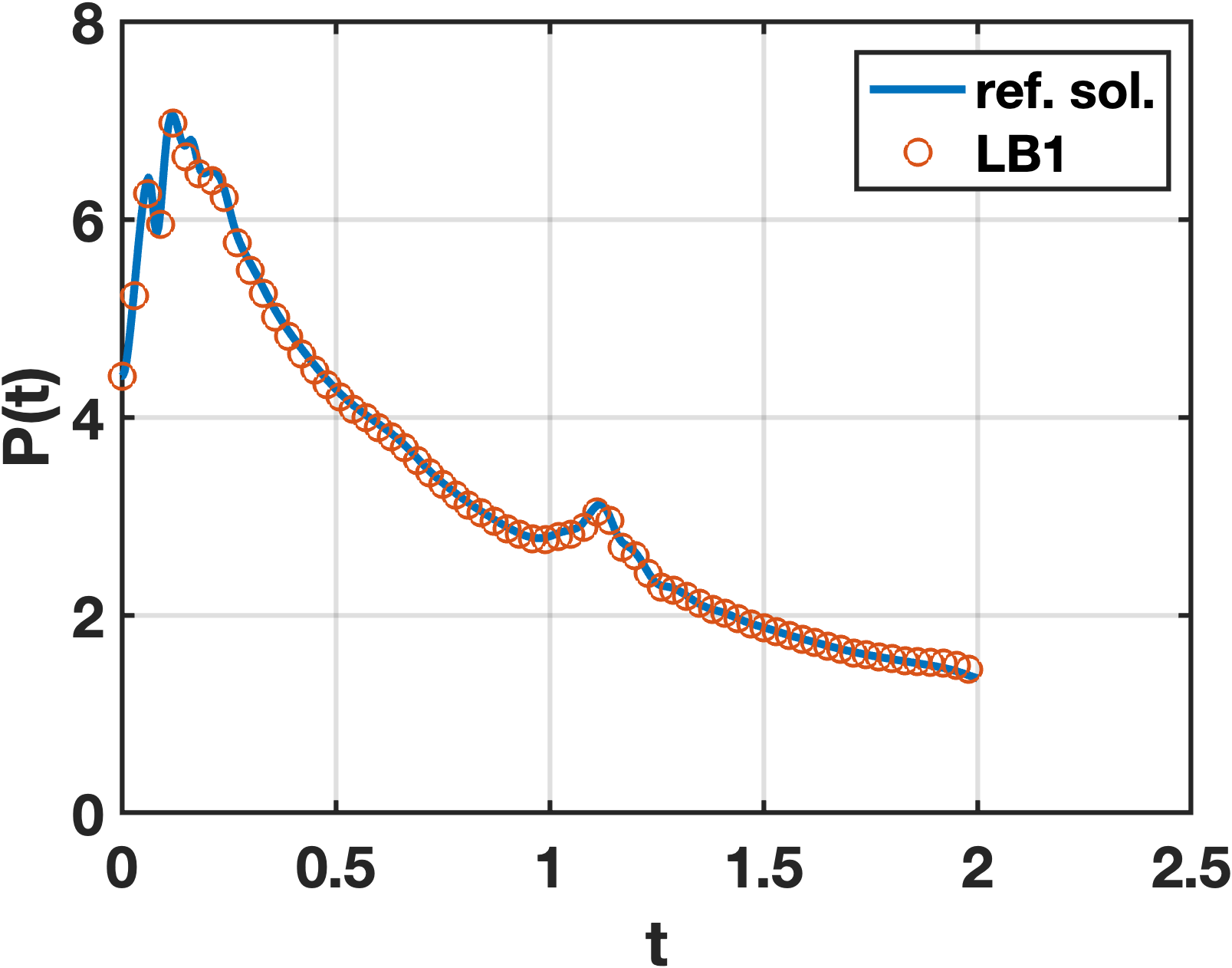}
        \\ \small (c) Palinstrophy $P(t)$
    \end{minipage}

    % Légende globale descriptive pour votre rapport/article
    \caption{Comparison of global physical quantities between the pseudo-spectral reference solution and the lattice Boltzmann (LB1) simulation: (a) kinetic energy, (b) enstrophy, and (c) palinstrophy.}
    \label{fig:comparaison_grandeurs_globales}
\end{figure}

\begin{figure}[htbp]
    \centering
    % --- Première figure : Champ de vorticité LBM ---
    \begin{minipage}[b]{0.50\textwidth}
        \centering
        % Remplacez 1.8 par la valeur de votre time_label si nécessaire
        \includegraphics[width=\textwidth]{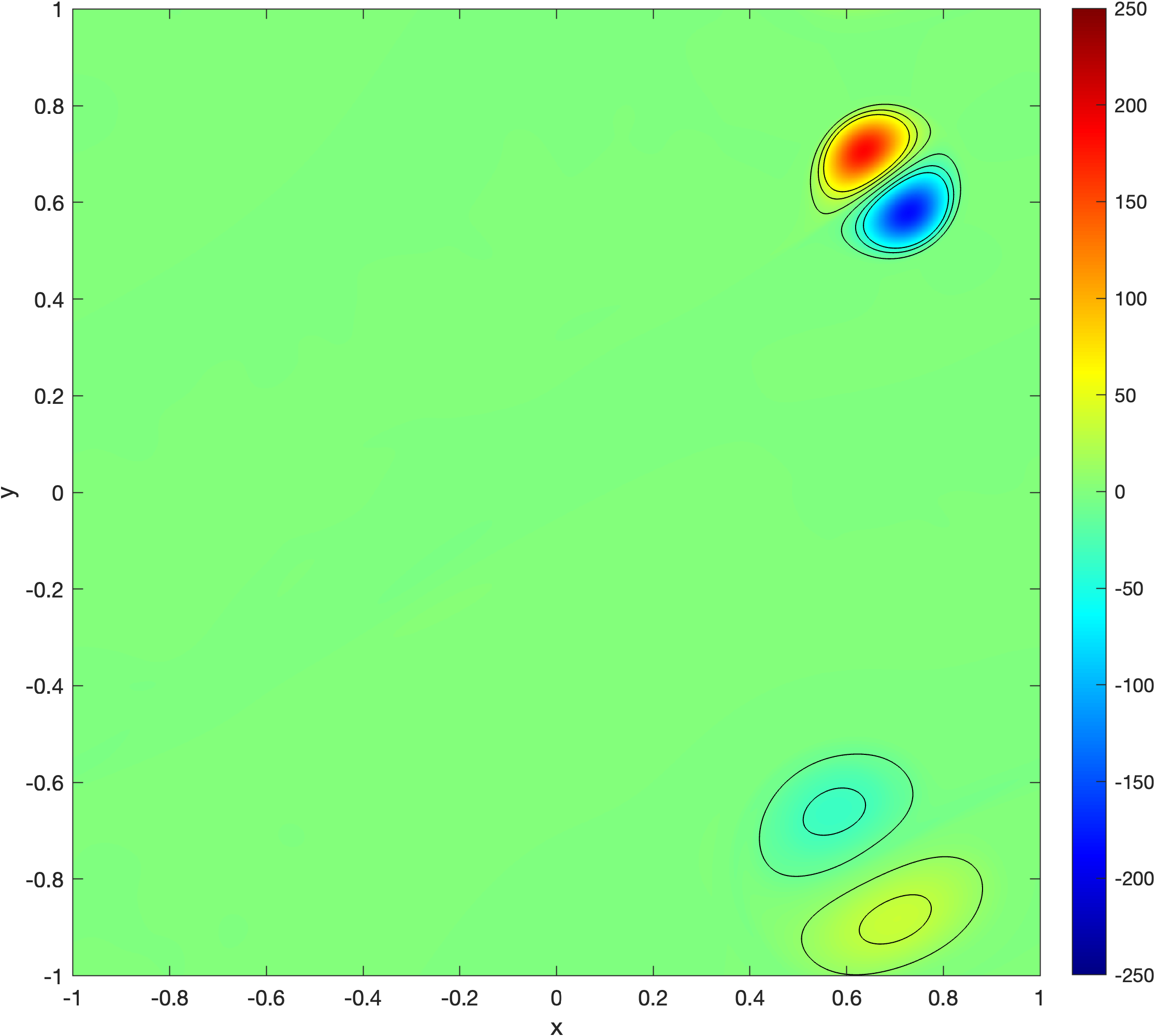}
        \\ \small (a) LB1 vorticity field
    \end{minipage}
    \hfill % Pousse les deux images aux extrémités pour un alignement propre
    % --- Deuxième figure : Contours de comparaison ---
    \begin{minipage}[b]{0.46\textwidth}
        \centering
        % Remplacez 1.8 par la valeur de votre time_label si nécessaire
        \includegraphics[width=\textwidth]{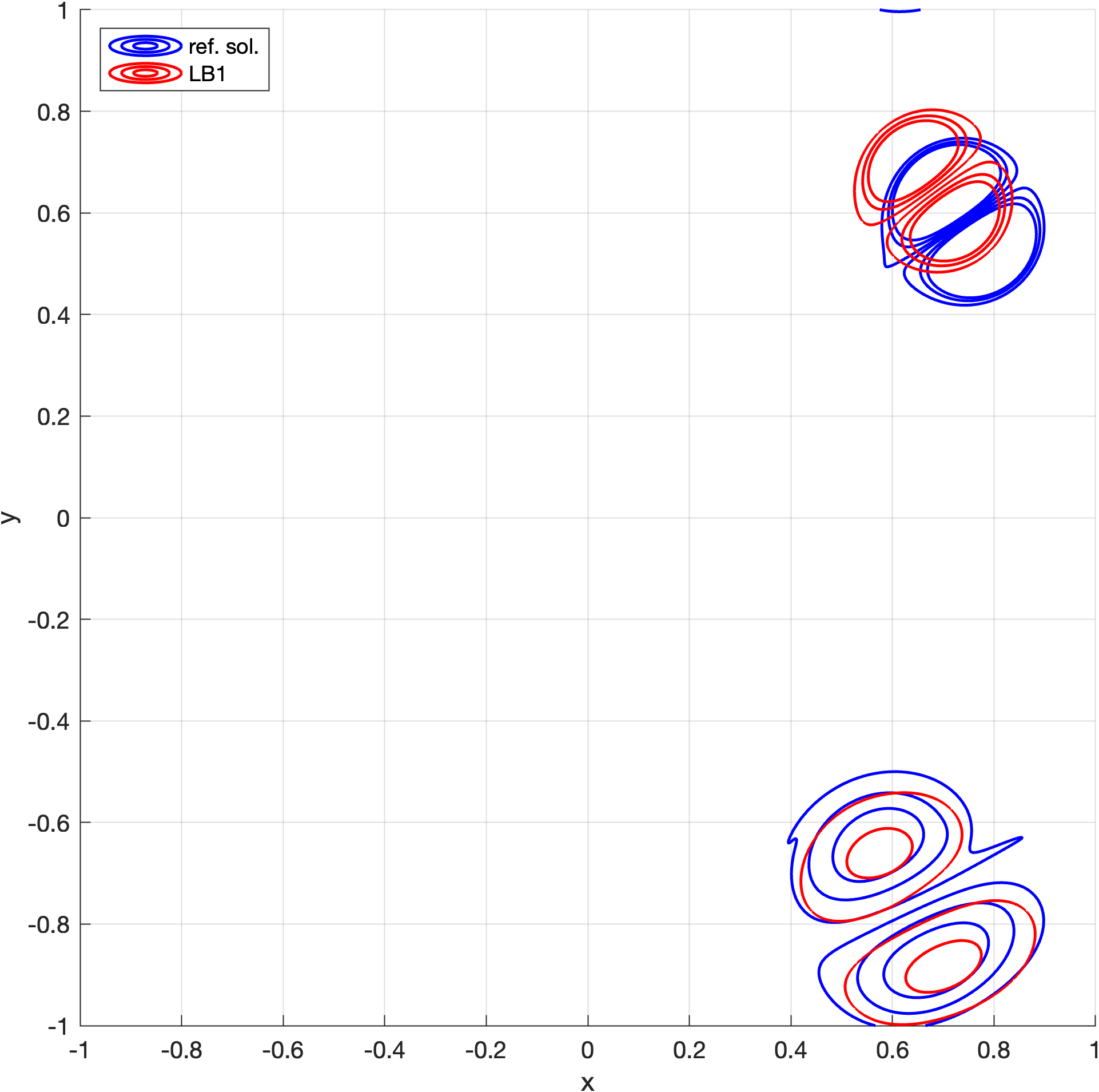}
        \\ \small (b) Spectral vs. LB1 contours
    \end{minipage}

    % Légende globale en anglais pour votre article
    \caption{Vorticity field results at $t = 1.8$: (a) computed vorticity field using the LB1 scheme, and (b) comparison of the vorticity contours between the pseudo-spectral reference solution (blue) and the LB1 scheme (red).}
    \label{fig:comparaison_vorticite_instantannee}
\end{figure}

\noindent First, we consider the LB1 scheme, which corresponds to the standard MRT model where the equilibrium distribution function is given by Eq.~(\ref{equi_qian}).
The Reynolds number $Re = 2500$ determines the relaxation rate $s_\nu$, while the other relaxation rates are chosen as follows: $s_e = 1.9$, $s_q = 1.93$, and $s_h = 1.94$.
Figure~\ref{fig:comparaison_grandeurs_globales} compares the kinetic energy, enstrophy, and palinstrophy computed by the LB1 scheme with those of the reference solution.
The results demonstrate a good agreement between the LB1 simulation and the reference solution.
Figure~\ref{fig:comparaison_vorticite_instantannee}(a) displays the vorticity field computed by the standard MRT scheme (LB1) at time $t=1.8$.
Figure~\ref{fig:comparaison_vorticite_instantannee}(b) compares the vorticity contours of the LB1 solution with those of the reference solution. It can be observed that the positions of the two vortices obtained with LB1 are shifted to the left compared to the reference solution.
This clearly highlights the lack of isotropy in the LB1 scheme, which is caused by the spurious terms $\mathcal{S}_x$ and $\mathcal{S}_y$ (see Eqs.~(\ref{S_X}) and~(\ref{S_Y})) that break Galilean invariance.
These results confirm those obtained by Lallemand \emph{et al.} in \cite{LDL24}.

%%%%%%%%%
%%%%%%%%%
%%%%%%%%%.      LB2
%%%%%%%%%
%%%%%%%%%

\subsection*{LB2 scheme}
\begin{figure}[htbp]
    \centering
    % --- Première figure : Comparaison Énergie ---
    \begin{minipage}[b]{0.32\textwidth}
        \centering
        \includegraphics[width=\textwidth]{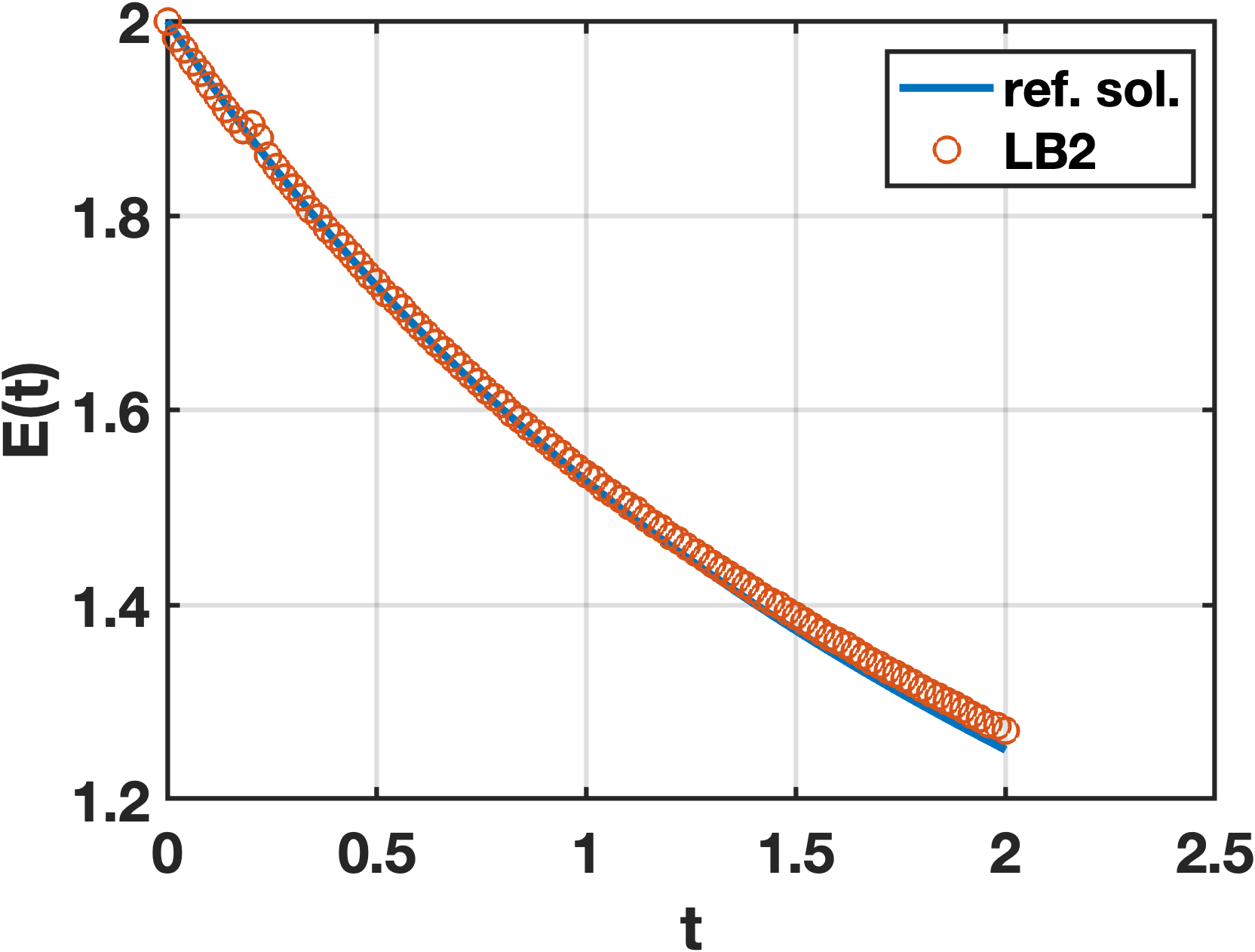}
        \\ \small (a) Kinetic energy $E(t)$
    \end{minipage}
    \hfill % Ajuste l'espace horizontal de manière homogène
    % --- Deuxième figure : Comparaison Enstrophie ---
    \begin{minipage}[b]{0.32\textwidth}
        \centering
        \includegraphics[width=\textwidth]{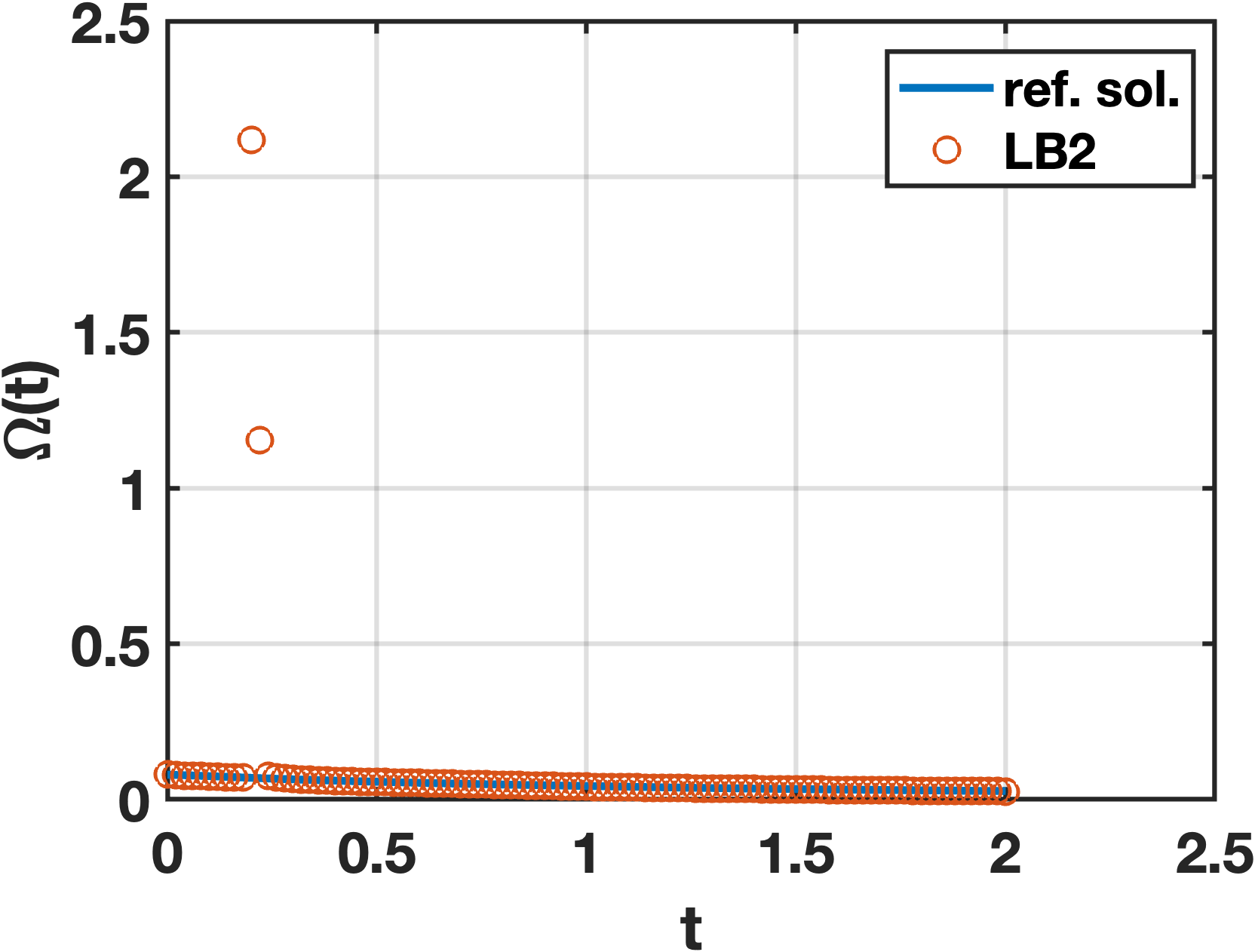}
        \\ \small (b) Enstrophy $\Omega(t)$
    \end{minipage}
    \hfill % Ajuste l'espace horizontal de manière homogène
    % --- Troisième figure : Comparaison Palinstrophie ---
    \begin{minipage}[b]{0.32\textwidth}
        \centering
        \includegraphics[width=\textwidth]{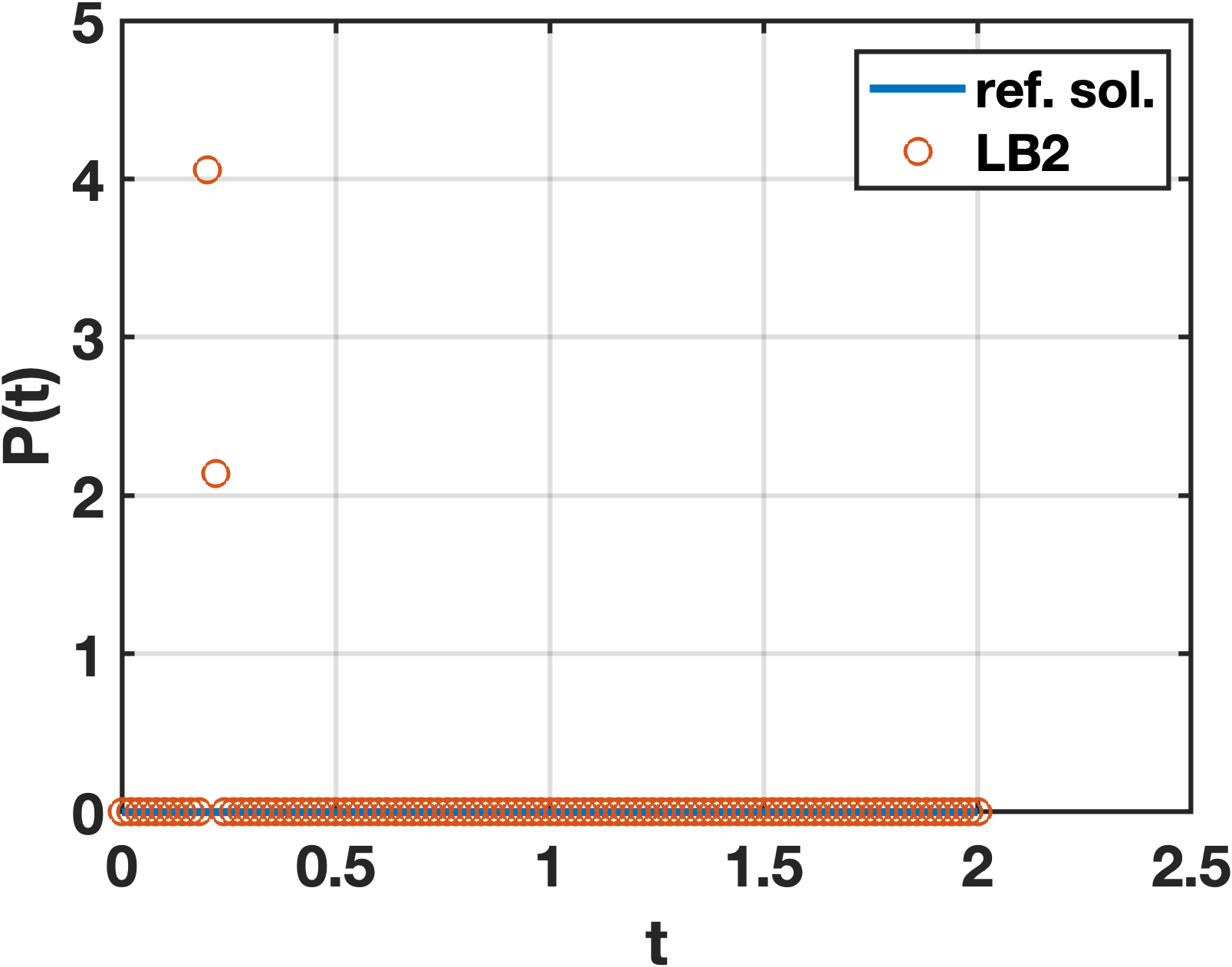}
        \\ \small (c) Palinstrophy $P(t)$
    \end{minipage}

    % Légende globale descriptive pour votre rapport/article
    \caption{Comparison of global physical quantities between the pseudo-spectral reference solution and the lattice Boltzmann (LB2) simulation: (a) kinetic energy, (b) enstrophy, and (c) palinstrophy.}
    \label{fig:comp_LB2}
\end{figure}

\noindent The LB2 scheme is identical to the above LB1 scheme, except for the equilibrium of the fourth-order moment $h$. Specifically, the non-linear terms of $h^{eq}$ are eliminated, yielding $h^{eq}=0$.
Given that $h^{eq}$ has no influence on the second-order equivalent PDEs, this choice is intended to isolate and study its specific effect on the oblique dipole benchmark.
Regarding the choice of relaxation rates, they are kept identical to those of the LB1 scheme.
Figure~\ref{fig:comp_LB2} compares the kinetic energy, enstrophy, and palinstrophy computed by the LB2 scheme with those of the reference solution.
The results demonstrate that the LB2 solution exhibits an unphysical behavior, with global quantities that significantly deviate from the reference solution.

\begin{figure}[htbp]
    \centering
    % --- Première figure : Champ de vorticité LBM ---
    \begin{minipage}[b]{0.50\textwidth}
        \centering
        % Remplacez 1.8 par la valeur de votre time_label si nécessaire
        \includegraphics[width=\textwidth]{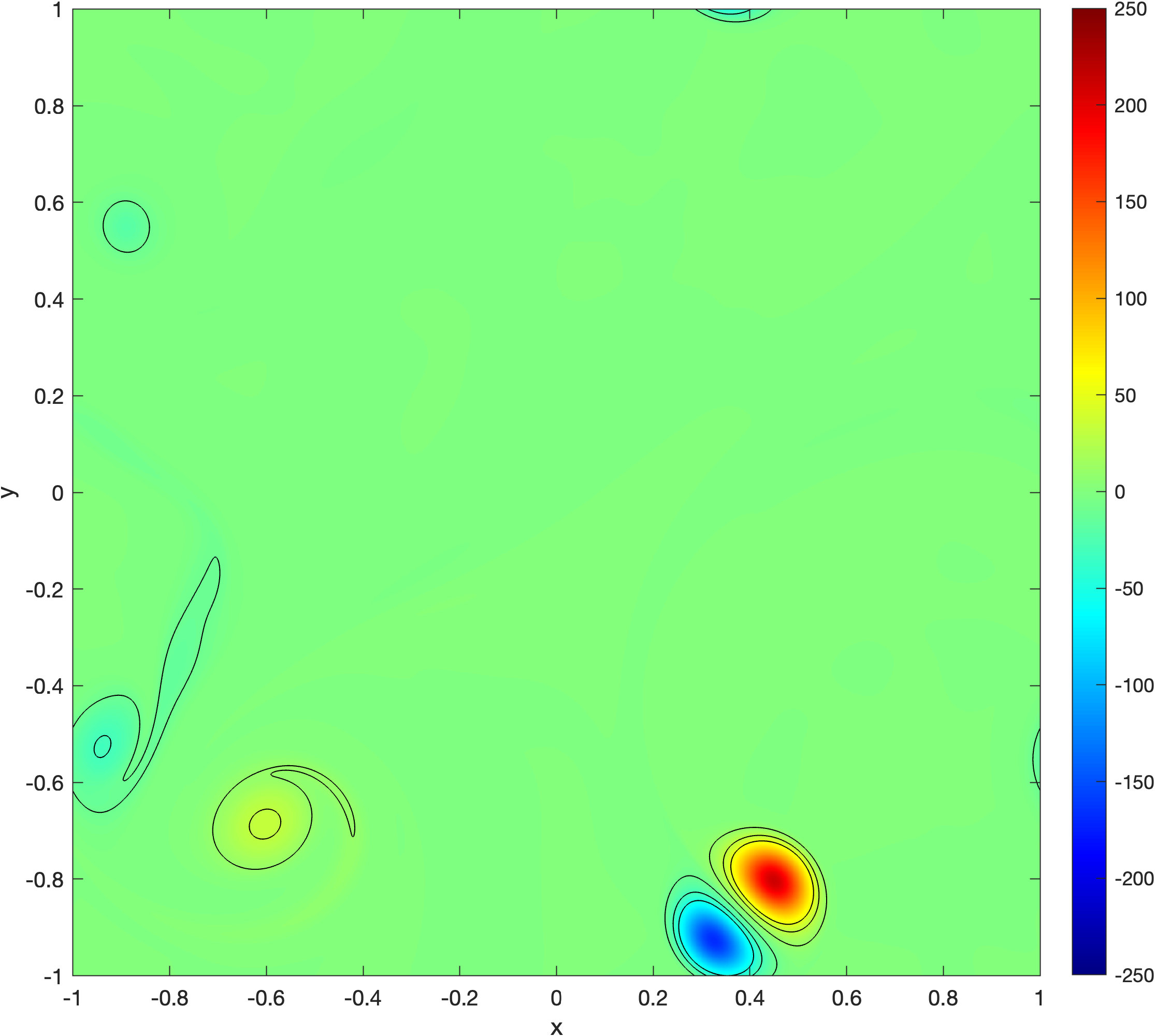}
        \\ \small (a) LB2 vorticity field
    \end{minipage}
    \hfill % Pousse les deux images aux extrémités pour un alignement propre
    % --- Deuxième figure : Contours de comparaison ---
    \begin{minipage}[b]{0.45\textwidth}
        \centering
        % Remplacez 1.8 par la valeur de votre time_label si nécessaire
        \includegraphics[width=\textwidth]{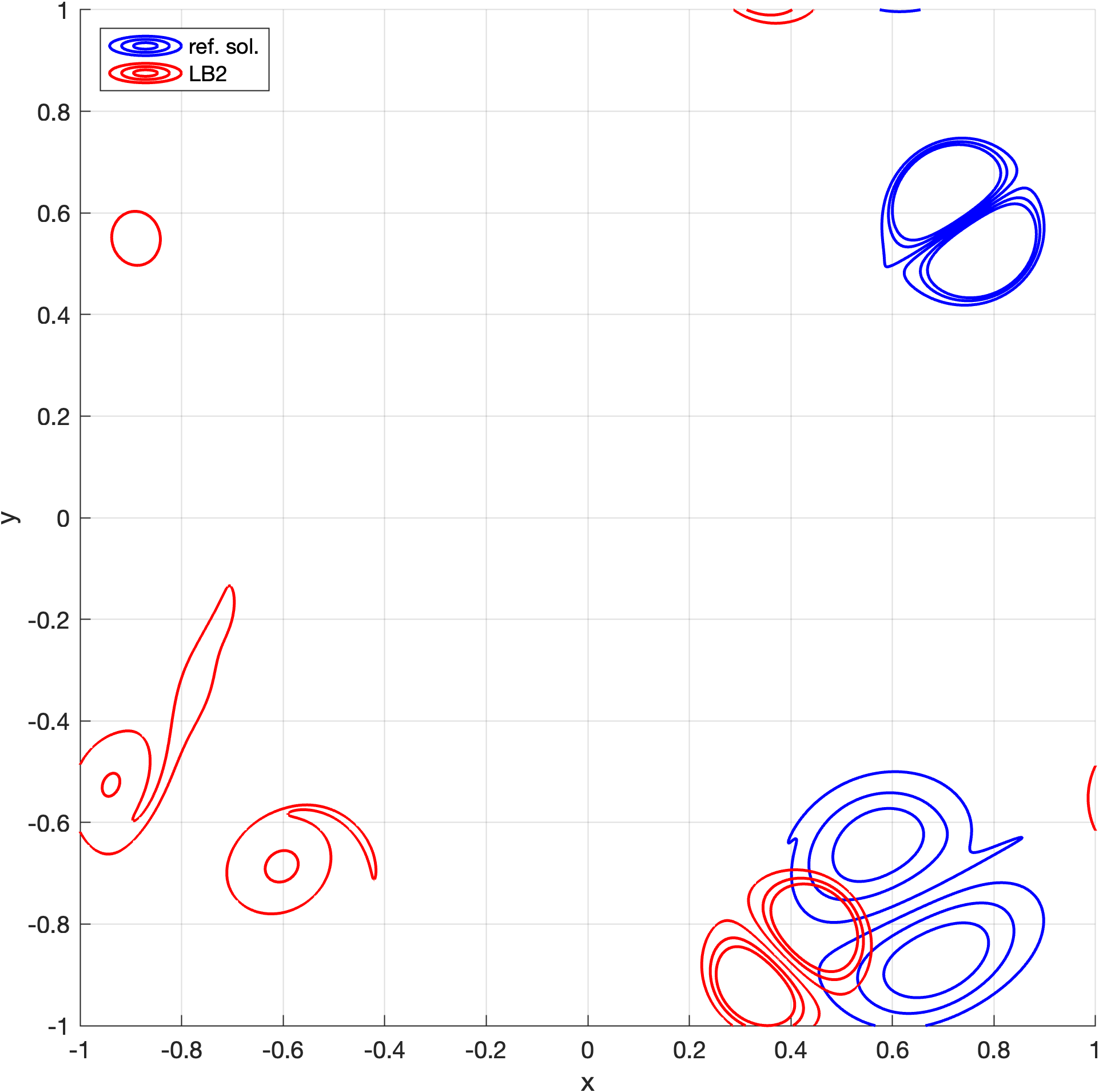}
        \\ \small (b) Spectral vs. LB2 contours
    \end{minipage}

    % Légende globale en anglais pour votre article
    \caption{Vorticity field results at $t = 1.8$: (a) computed vorticity field using the LB2 scheme, and (b) comparison of the vorticity contours between the pseudo-spectral reference solution (blue) and the LB2 scheme (red).}
    \label{fig:comp_vor_LB2}
\end{figure}

Figure~\ref{fig:comp_vor_LB2}(a) displays the vorticity field computed by the scheme (LB2) at time $t=1.8$.
Figure~\ref{fig:comp_vor_LB2}(b) compares the vorticity contours of the LB2 solution with those of the reference solution.
The figure clearly shows that the coherent structure of the two vortices is no longer preserved, leading to a significant deviation from the reference solution.
This demonstrates that the non-linear terms in the definition of $h^{eq}$ are fundamental, even though they only affect the equivalent partial differential equations (PDEs) at higher orders.
It is worth noting that the choice of $h^{eq}$ has no noticeable effect on canonical benchmarks such as the Poiseuille flow, the lid-driven cavity, or the Taylor--Green vortex.
Consequently, the present benchmark confirms that the full equilibrium distribution function described by Eq.~(\ref{equi_qian}), originally introduced in~\cite{QHL92},
is required not only for second-order consistency but also for its appropriate higher-order contributions.

%%%%%%%%%
%%%%%%%%%
%%%%%%%%%.      LB3
%%%%%%%%%
%%%%%%%%%

\subsection*{LB3 scheme}
\begin{figure}[htbp]
    \centering
    % --- Première figure : Comparaison Énergie ---
    \begin{minipage}[b]{0.32\textwidth}
        \centering
        \includegraphics[width=\textwidth]{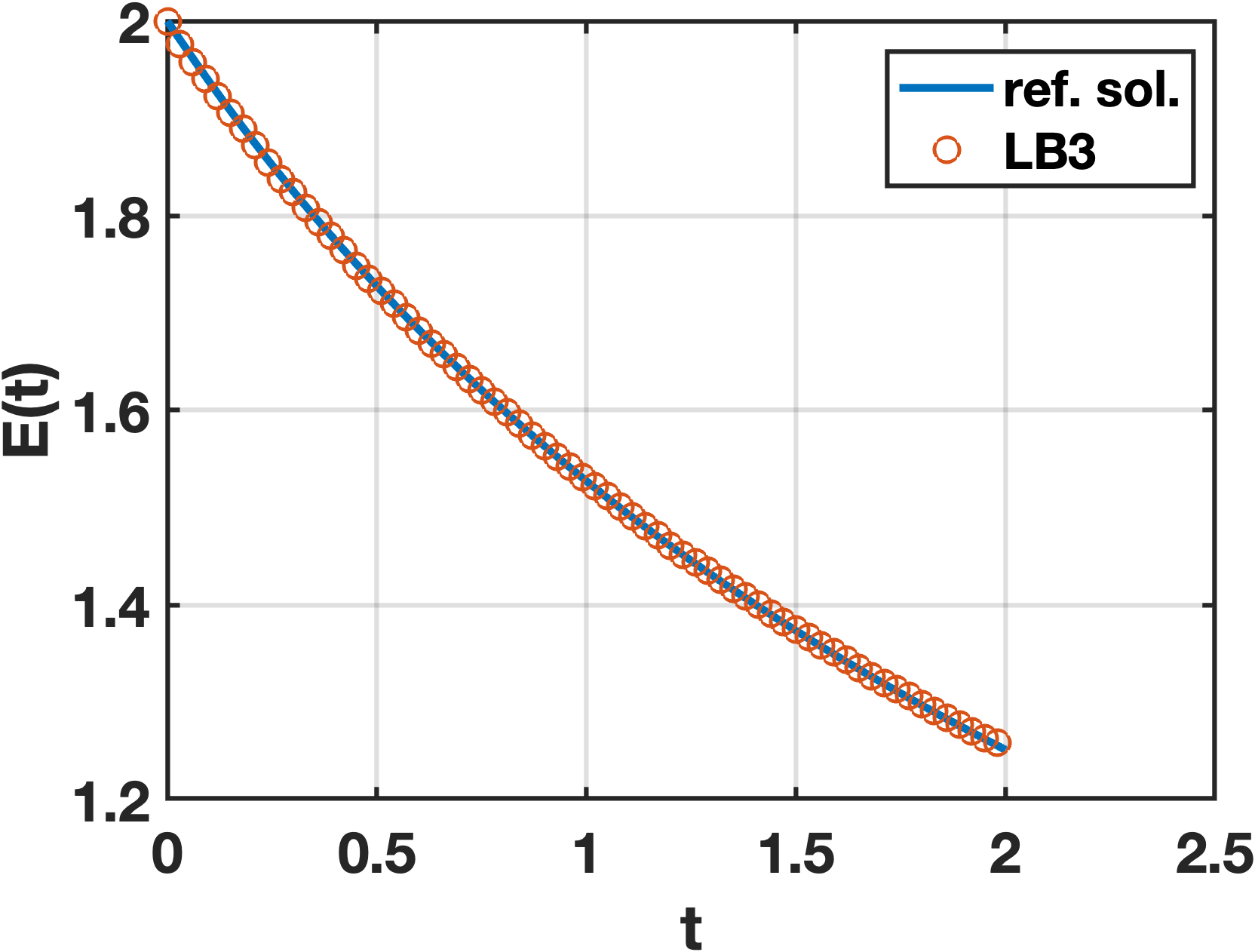}
        \\ \small (a) Kinetic energy $E(t)$
    \end{minipage}
    \hfill % Ajuste l'espace horizontal de manière homogène
    % --- Deuxième figure : Comparaison Enstrophie ---
    \begin{minipage}[b]{0.32\textwidth}
        \centering
        \includegraphics[width=\textwidth]{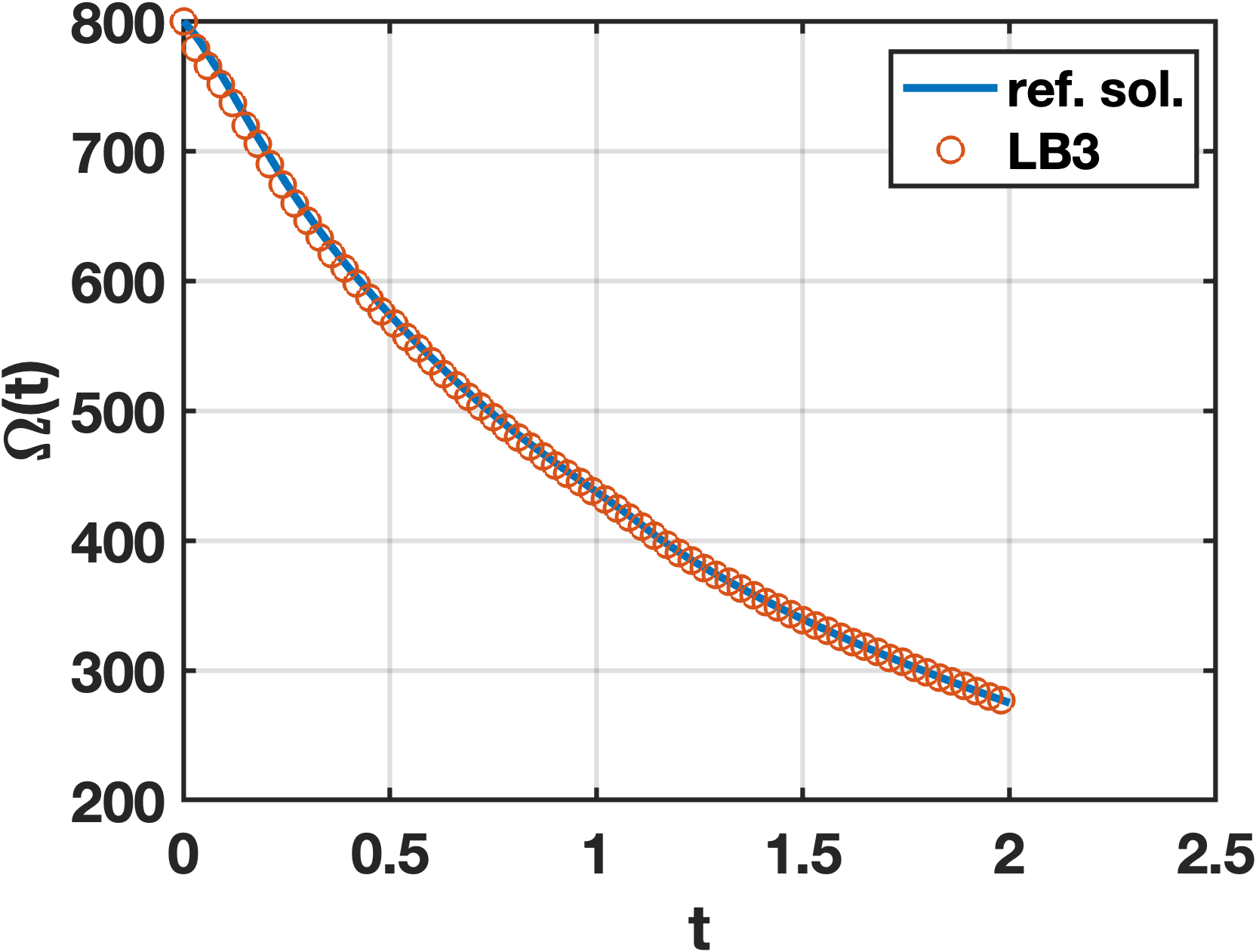}
        \\ \small (b) Enstrophy $\Omega(t)$
    \end{minipage}
    \hfill % Ajuste l'espace horizontal de manière homogène
    % --- Troisième figure : Comparaison Palinstrophie ---
    \begin{minipage}[b]{0.32\textwidth}
        \centering
        \includegraphics[width=\textwidth]{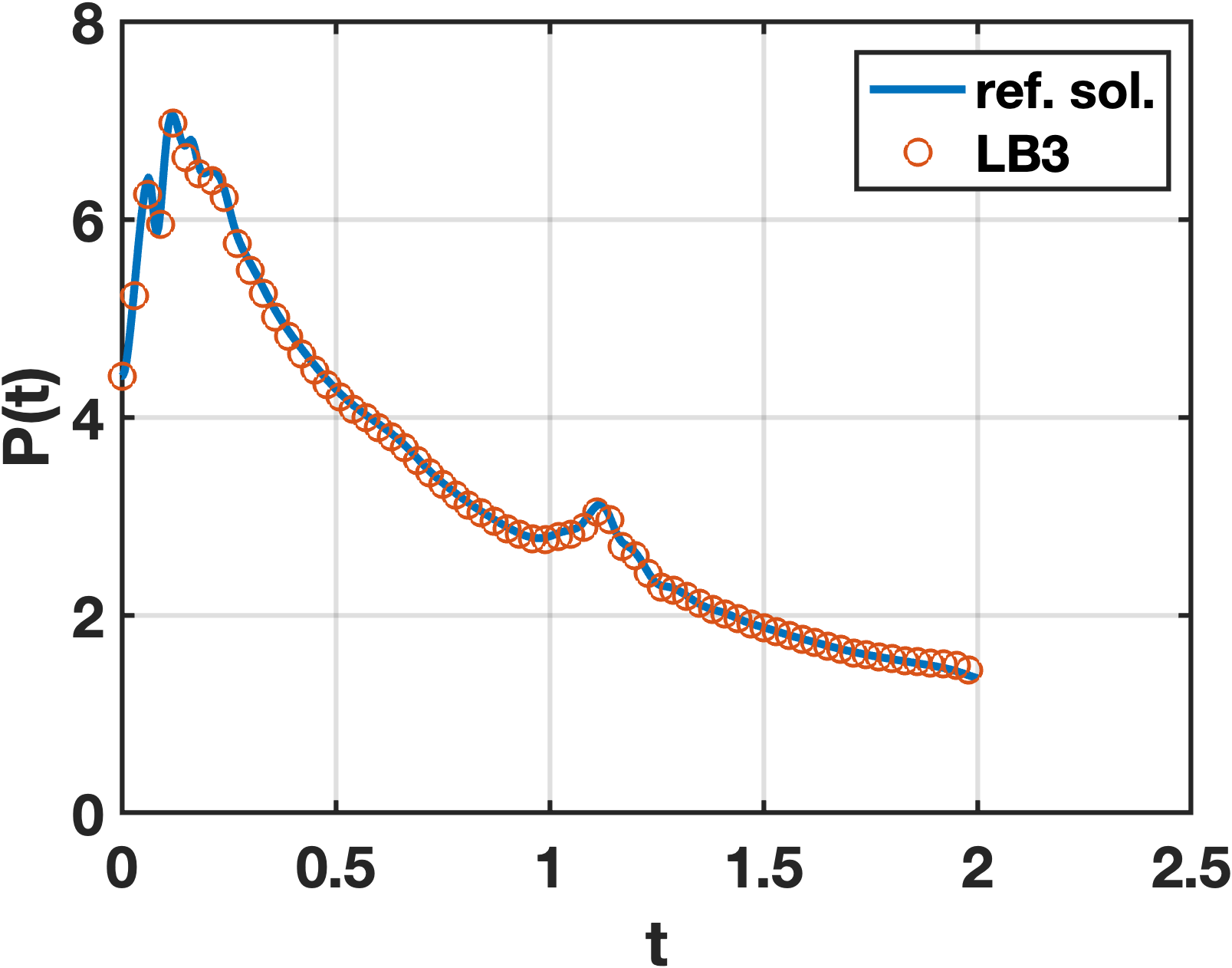}
        \\ \small (c) Palinstrophy $P(t)$
    \end{minipage}

    % Légende globale descriptive pour votre rapport/article
    \caption{Comparison of global physical quantities between the pseudo-spectral reference solution and the lattice Boltzmann (LB3) simulation: (a) kinetic energy, (b) enstrophy, and (c) palinstrophy.}
    \label{fig:comp_LB3}
\end{figure}

\begin{figure}[htbp]
    \centering
    % --- Première figure : Champ de vorticité LBM ---
    \begin{minipage}[b]{0.50\textwidth}
        \centering
        % Remplacez 1.8 par la valeur de votre time_label si nécessaire
        \includegraphics[width=\textwidth]{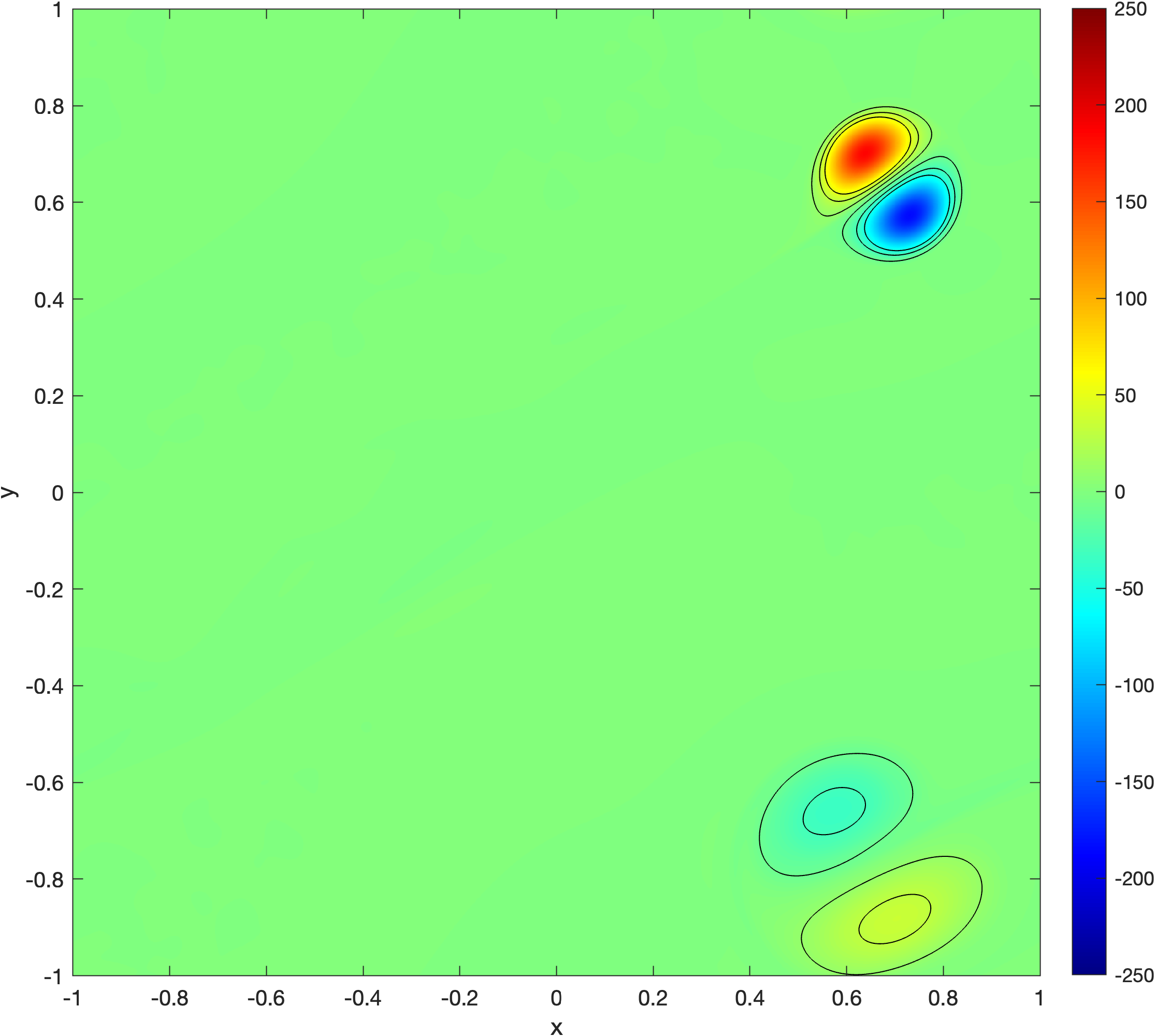}
        \\ \small (a) LB3 vorticity field
    \end{minipage}
    \hfill % Pousse les deux images aux extrémités pour un alignement propre
    % --- Deuxième figure : Contours de comparaison ---
    \begin{minipage}[b]{0.45\textwidth}
        \centering
        % Remplacez 1.8 par la valeur de votre time_label si nécessaire
        \includegraphics[width=\textwidth]{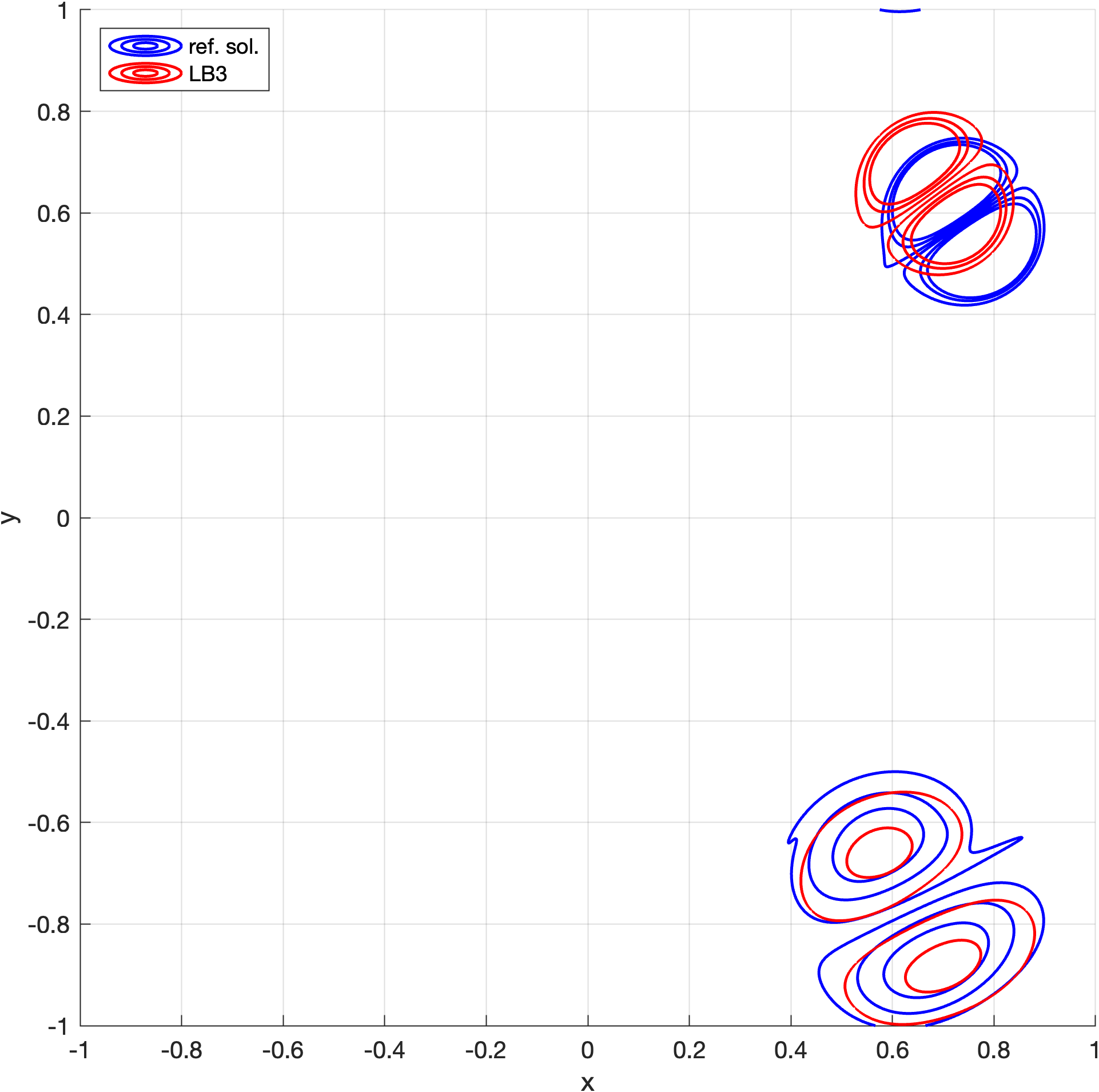}
        \\ \small (b) Spectral vs. LB3 contours
    \end{minipage}

    % Légende globale en anglais pour votre article
    \caption{Vorticity field results at $t = 1.8$: (a) computed vorticity field using the LB3 scheme, and (b) comparison of the vorticity contours between the pseudo-spectral reference solution (blue) and the LB3 scheme (red).}
    \label{fig:comp_vor_LB3}
\end{figure}

\noindent  We consider the LB3 scheme, which corresponds to the projected LBM scheme.
Again, the Reynolds number $Re = 2500$ determines the relaxation rate $s_x$, while the only remaining free relaxation rate is chosen as $s_e = 1.99$.
It is worth noting that this scheme allows for a smaller bulk viscosity $\zeta$ compared to the standard MRT (LB1) scheme.
This is due to the enhanced stability of the projected formulation; indeed, for the considered value of $s_e = 1.99$, the LB1 scheme diverges.
Figure~\ref{fig:comp_LB3} compares the kinetic energy, enstrophy, and palinstrophy computed by the LB3 scheme with those of the reference solution.
%The results demonstrate a good agreement between the LB1 simulation and the reference solution.
Figure~\ref{fig:comp_vor_LB3}(a) displays the vorticity field computed by the scheme LB3 at time $t=1.8$.
Figure~\ref{fig:comp_vor_LB3}(b) compares the vorticity contours of the LB2 solution with those of the reference solution.
As the LB1 scheme, it can be observed that the positions of the two vortices obtained with LB3 are shifted to the left compared to the reference solution.

%\newpage

%%%%%%%%%
%%%%%%%%%
%%%%%%%%%.      LB4
%%%%%%%%%
%%%%%%%%%

\subsection*{LB4 scheme}
\begin{figure}[htbp]
    \centering
    % --- Première figure : Comparaison Énergie ---
    \begin{minipage}[b]{0.32\textwidth}
        \centering
        \includegraphics[width=\textwidth]{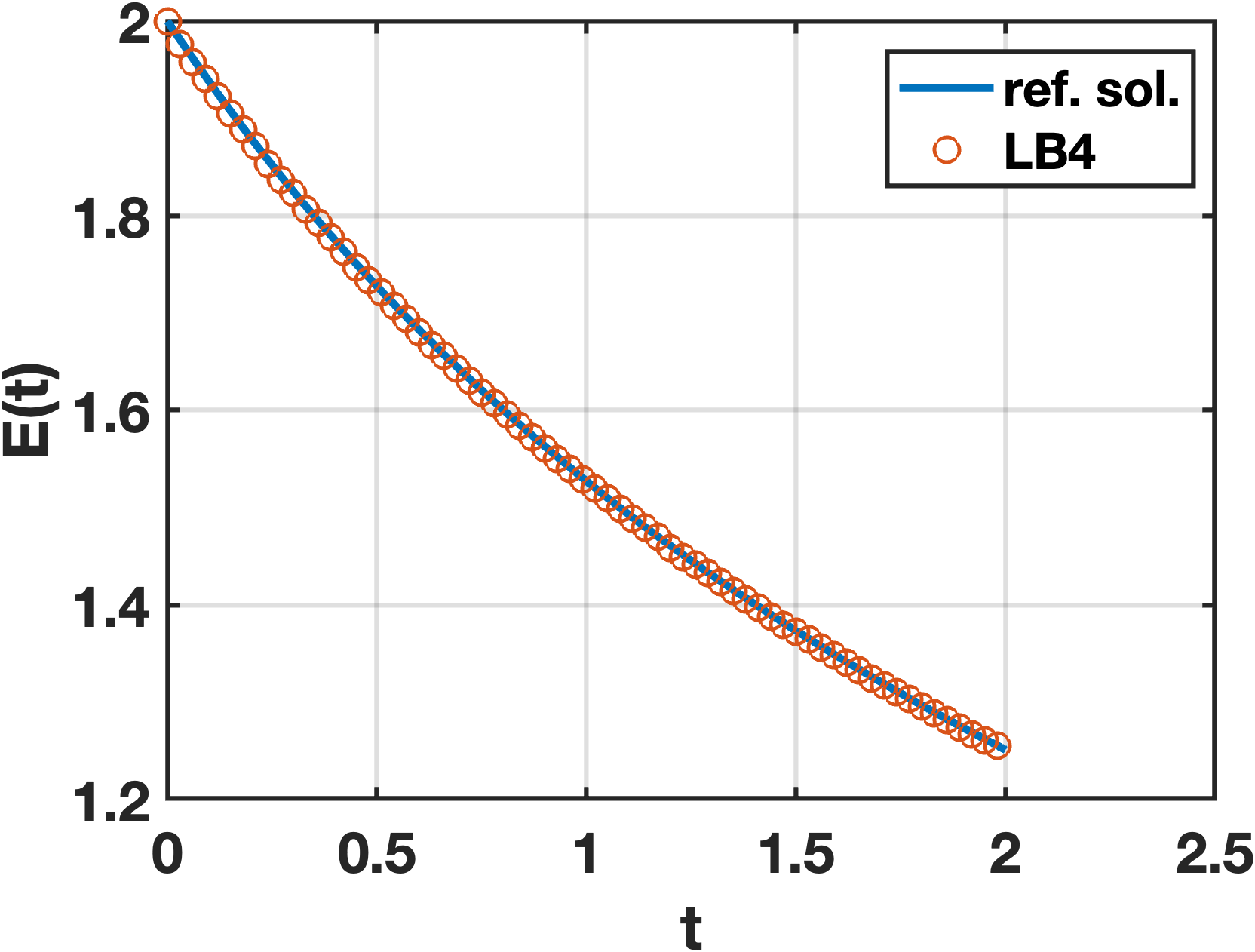}
        \\ \small (a) Kinetic energy $E(t)$
    \end{minipage}
    \hfill % Ajuste l'espace horizontal de manière homogène
    % --- Deuxième figure : Comparaison Enstrophie ---
    \begin{minipage}[b]{0.32\textwidth}
        \centering
        \includegraphics[width=\textwidth]{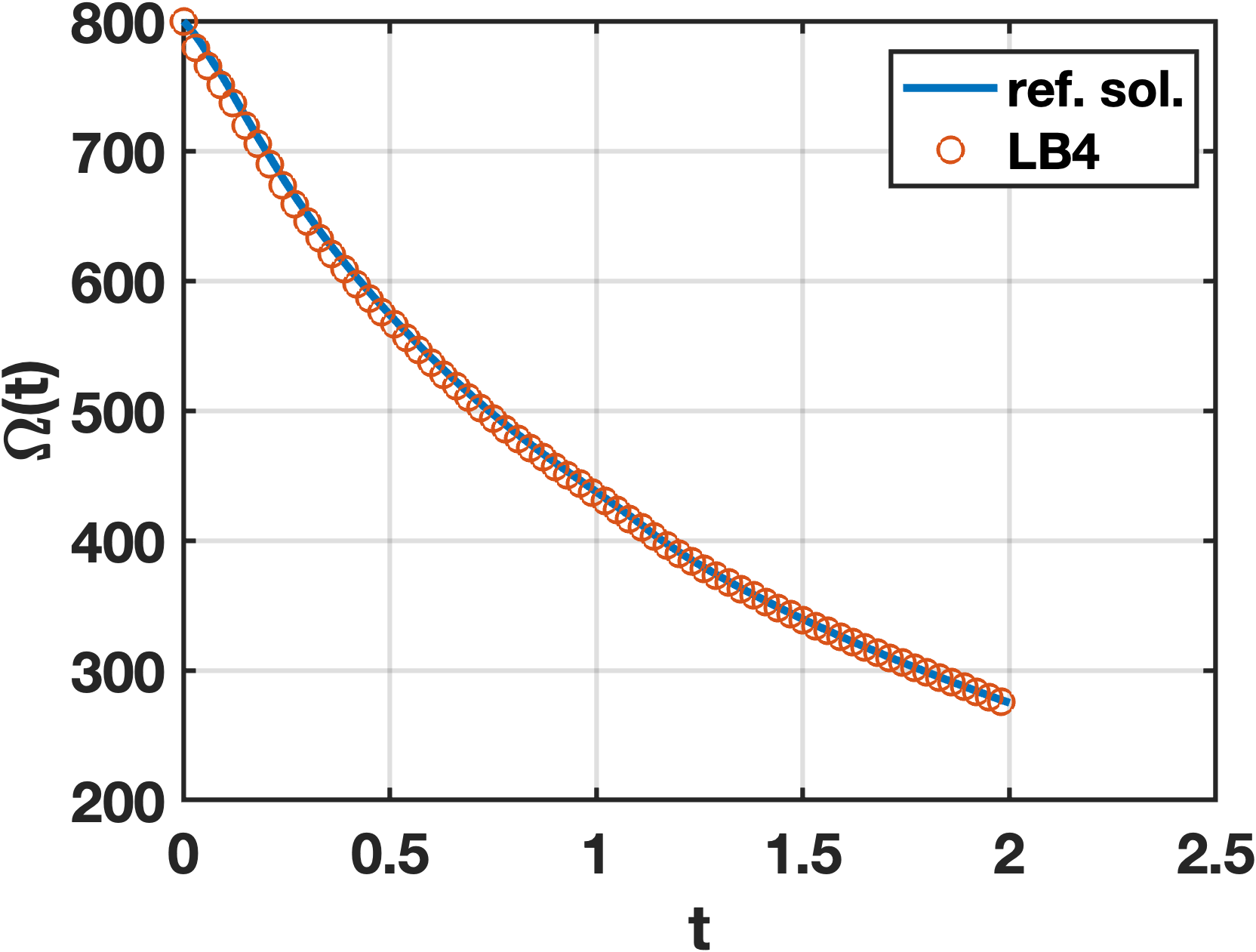}
        \\ \small (b) Enstrophy $\Omega(t)$
    \end{minipage}
    \hfill % Ajuste l'espace horizontal de manière homogène
    % --- Troisième figure : Comparaison Palinstrophie ---
    \begin{minipage}[b]{0.32\textwidth}
        \centering
        \includegraphics[width=\textwidth]{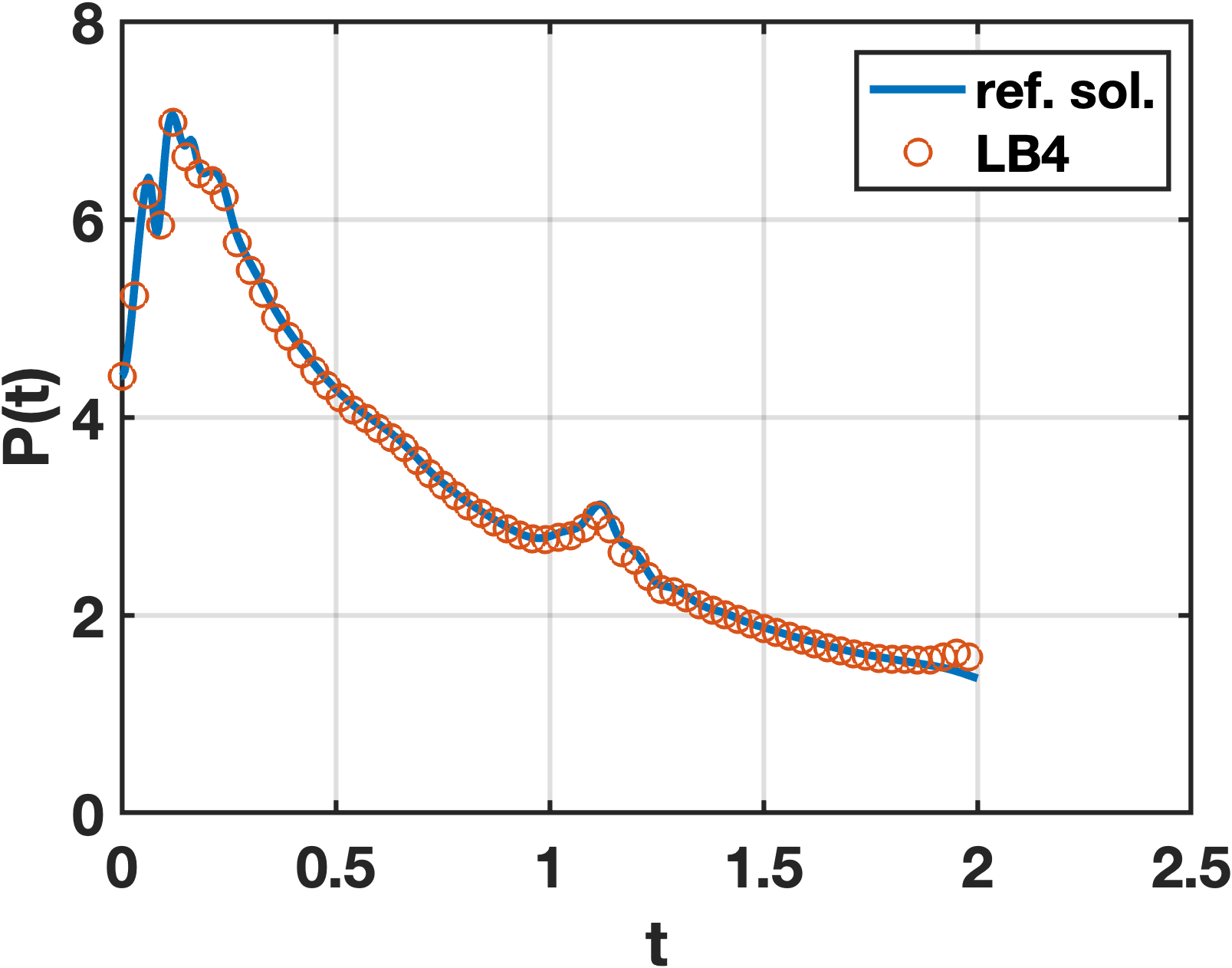}
        \\ \small (c) Palinstrophy $P(t)$
    \end{minipage}

    % Légende globale descriptive pour votre rapport/article
    \caption{Comparison of global physical quantities between the pseudo-spectral reference solution and the lattice Boltzmann (LB4) simulation: (a) kinetic energy, (b) enstrophy, and (c) palinstrophy.}
    \label{fig:comp_LB4}
\end{figure}

\begin{figure}[htbp]
    \centering
    % --- Première figure : Champ de vorticité LBM ---
    \begin{minipage}[b]{0.50\textwidth}
        \centering
        % Remplacez 1.8 par la valeur de votre time_label si nécessaire
        \includegraphics[width=\textwidth]{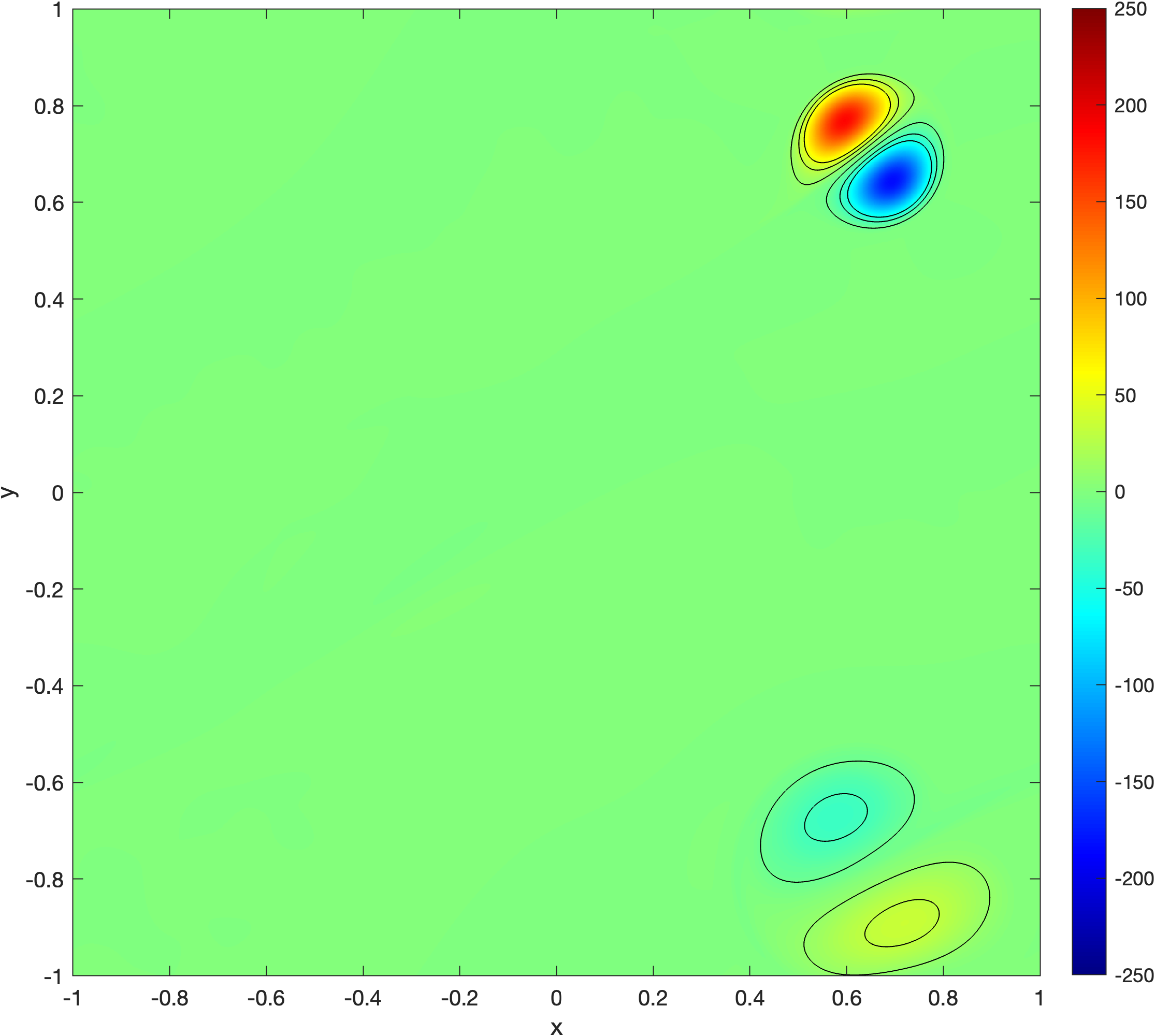}
        \\ \small (a) LB4 vorticity field
    \end{minipage}
    \hfill % Pousse les deux images aux extrémités pour un alignement propre
    % --- Deuxième figure : Contours de comparaison ---
    \begin{minipage}[b]{0.45\textwidth}
        \centering
        % Remplacez 1.8 par la valeur de votre time_label si nécessaire
        \includegraphics[width=\textwidth]{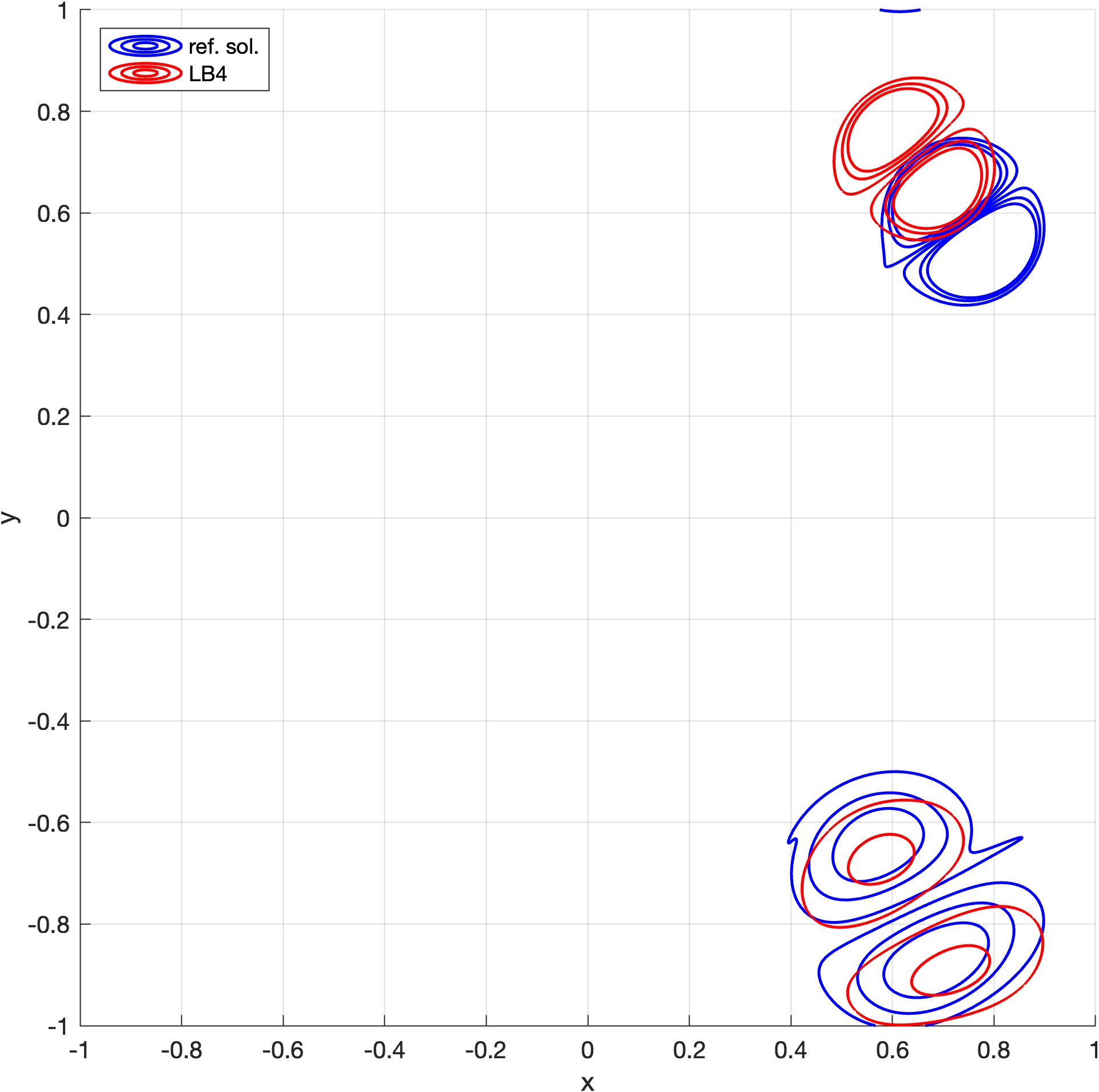}
        \\ \small (b) Spectral vs. LB4 contours
    \end{minipage}

    % Légende globale en anglais pour votre article
    \caption{Vorticity field results at $t = 1.8$: (a) computed vorticity field using the LB4 scheme, and (b) comparison of the vorticity contours between the pseudo-spectral reference solution (blue) and the LB4 scheme (red).}
    \label{fig:comp_vor_LB4}
\end{figure}

\noindent The LB4 scheme is identical to the LB1 scheme, except that it incorporates the Dubois equilibrium for $q_x$ and $q_y$, as given in Eq.~(\ref{du_eq}).
The relaxation rate $s_x$ is fixed to achieve a Reynolds number of $Re=2500$, while the remaining rates are set to $s_e=1.9$, $s_q=1.98$, and $s_h=1.9$.
Figure~\ref{fig:comp_LB4} compares the kinetic energy, enstrophy, and palinstrophy computed by the LB4 scheme with those of the reference solution.
%The results demonstrate a good agreement between the LB1 simulation and the reference solution.
Figure~\ref{fig:comp_vor_LB4}(a) displays the vorticity field computed by the scheme (LB4) at time $t=1.8$.
Figure~\ref{fig:comp_vor_LB4}(b) compares the vorticity contours of the LB3 solution with those of the reference solution.
%As the LB1 scheme, it can be observed that the positions of the two vortices obtained with LB2 are shifted to the left compared to the reference solution.
%On note ici que meme si les termes spurious $\mathcal{S}_x$ and $\mathcal{S}_y$ (see Eqs.~(\ref{S_X_FD}) and~(\ref{S_Y_FD})) ont mois de termes que que dans le cas du schema LB1 mais on vois qu'on a un schema non isotrope.
It is worth noting here that even though the spurious terms $\mathcal{S}_x$ and $\mathcal{S}_y$ (see Eqs.~(\ref{S_X_FD}) and~(\ref{S_Y_FD})) contain fewer terms than in the case of the LB1 scheme, the resulting simulation still exhibits a non-isotropic behavior.
Furthermore, the choice of the free relaxation rates $s_e$, $s_q$, and $s_h$ is primarily made to ensure the stability of the scheme.
Additionally, these parameters can be carefully selected to improve the isotropy of the model, which is precisely the strategy adopted for the subsequent LB5 scheme.

%%%%%%%%%
%%%%%%%%%
%%%%%%%%%.      LB5
%%%%%%%%%
%%%%%%%%%

\subsection*{LB5 scheme}
\begin{figure}[htbp]
    \centering
    % --- Première figure : Comparaison Énergie ---
    \begin{minipage}[b]{0.32\textwidth}
        \centering
        \includegraphics[width=\textwidth]{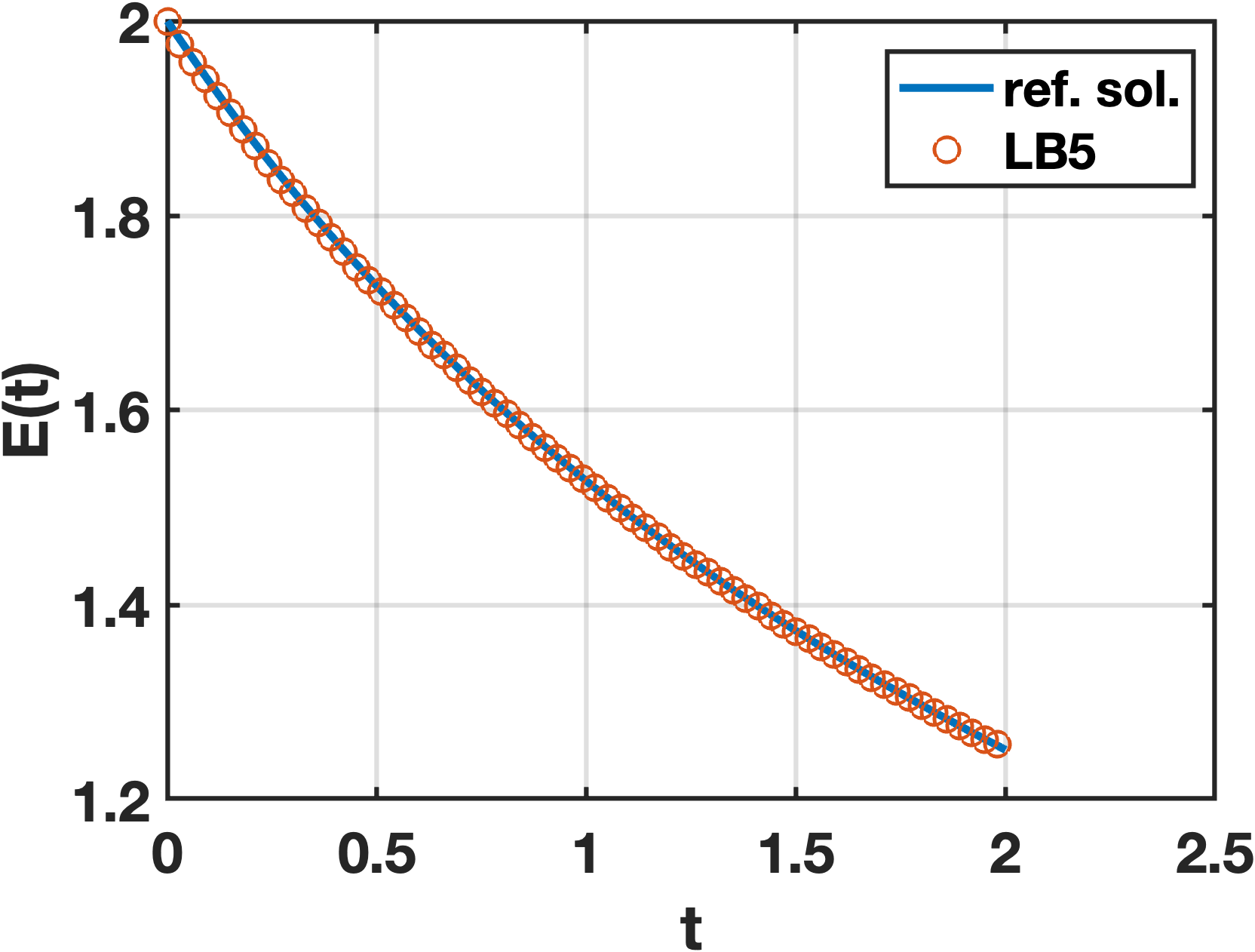}
        \\ \small (a) Kinetic energy $E(t)$
    \end{minipage}
    \hfill % Ajuste l'espace horizontal de manière homogène
    % --- Deuxième figure : Comparaison Enstrophie ---
    \begin{minipage}[b]{0.32\textwidth}
        \centering
        \includegraphics[width=\textwidth]{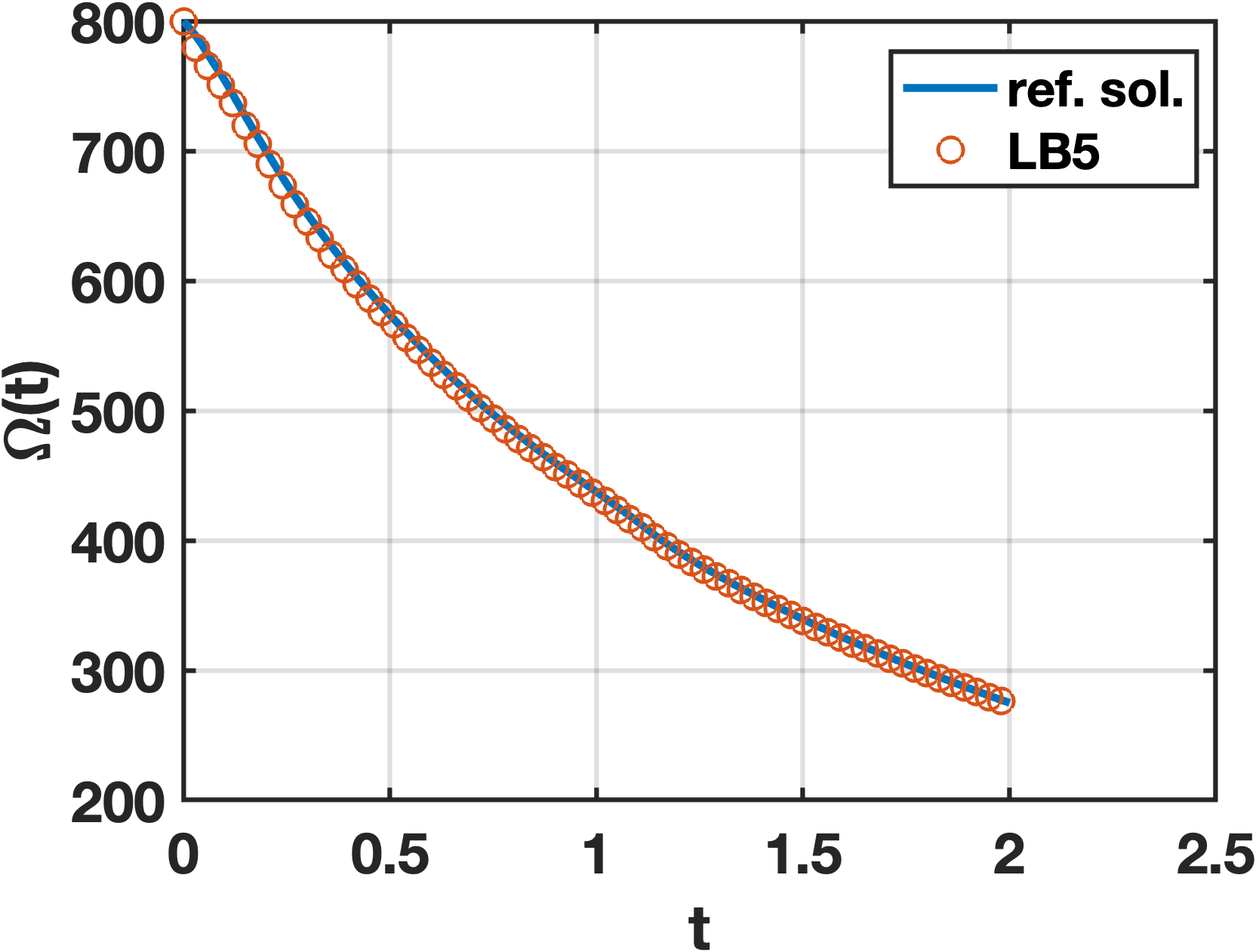}
        \\ \small (b) Enstrophy $\Omega(t)$
    \end{minipage}
    \hfill % Ajuste l'espace horizontal de manière homogène
    % --- Troisième figure : Comparaison Palinstrophie ---
    \begin{minipage}[b]{0.32\textwidth}
        \centering
        \includegraphics[width=\textwidth]{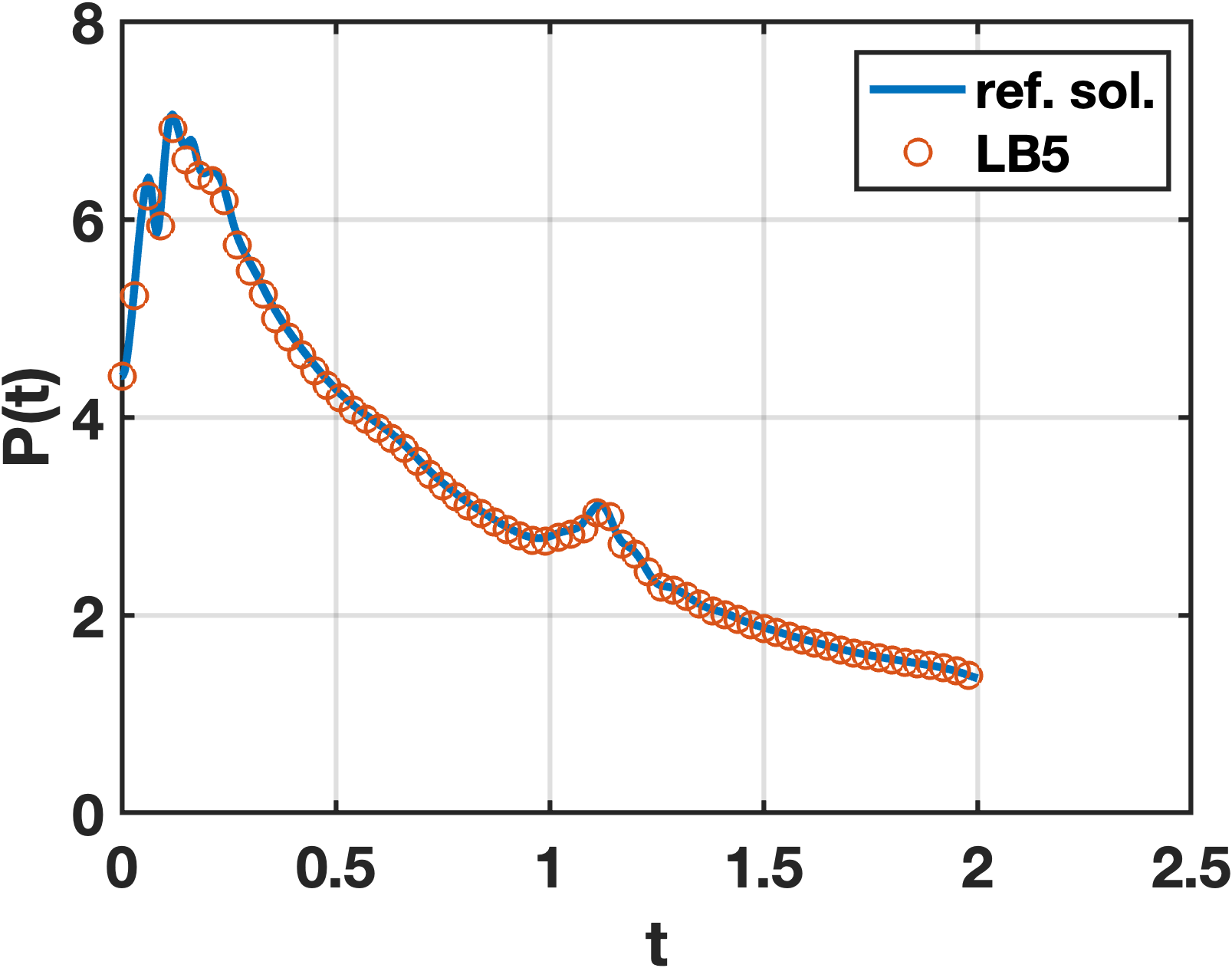}
        \\ \small (c) Palinstrophy $P(t)$
    \end{minipage}

    % Légende globale descriptive pour votre rapport/article
    \caption{Comparison of global physical quantities between the pseudo-spectral reference solution and the lattice Boltzmann (LB5) simulation: (a) kinetic energy, (b) enstrophy, and (c) palinstrophy.}
    \label{fig:comp_LB5}
\end{figure}

\begin{figure}[htbp]
    \centering
    % --- Première figure : Champ de vorticité LBM ---
    \begin{minipage}[b]{0.50\textwidth}
        \centering
        % Remplacez 1.8 par la valeur de votre time_label si nécessaire
        \includegraphics[width=\textwidth]{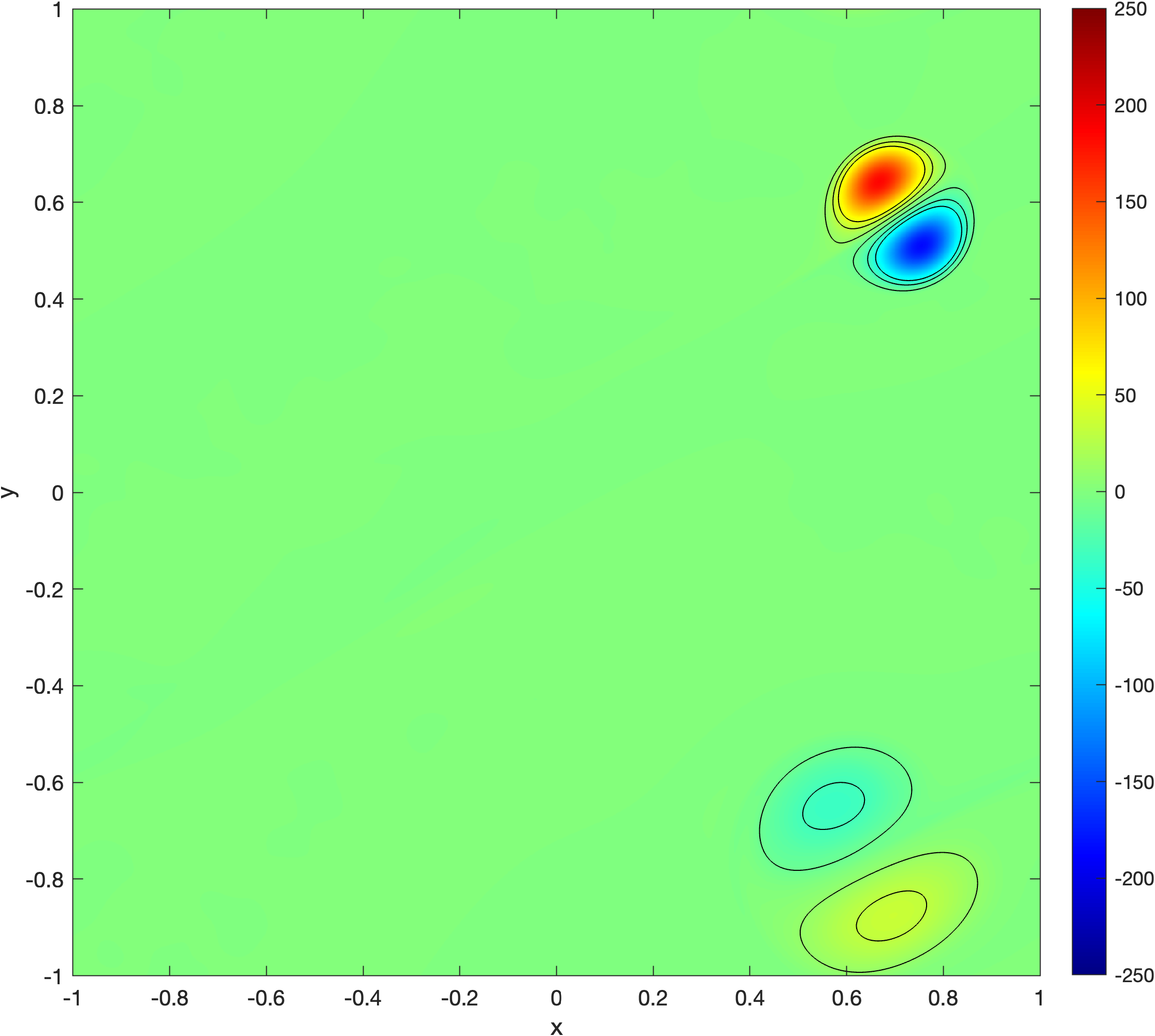}
        \\ \small (a) LB5 vorticity field
    \end{minipage}
    \hfill % Pousse les deux images aux extrémités pour un alignement propre
    % --- Deuxième figure : Contours de comparaison ---
    \begin{minipage}[b]{0.45\textwidth}
        \centering
        % Remplacez 1.8 par la valeur de votre time_label si nécessaire
        \includegraphics[width=\textwidth]{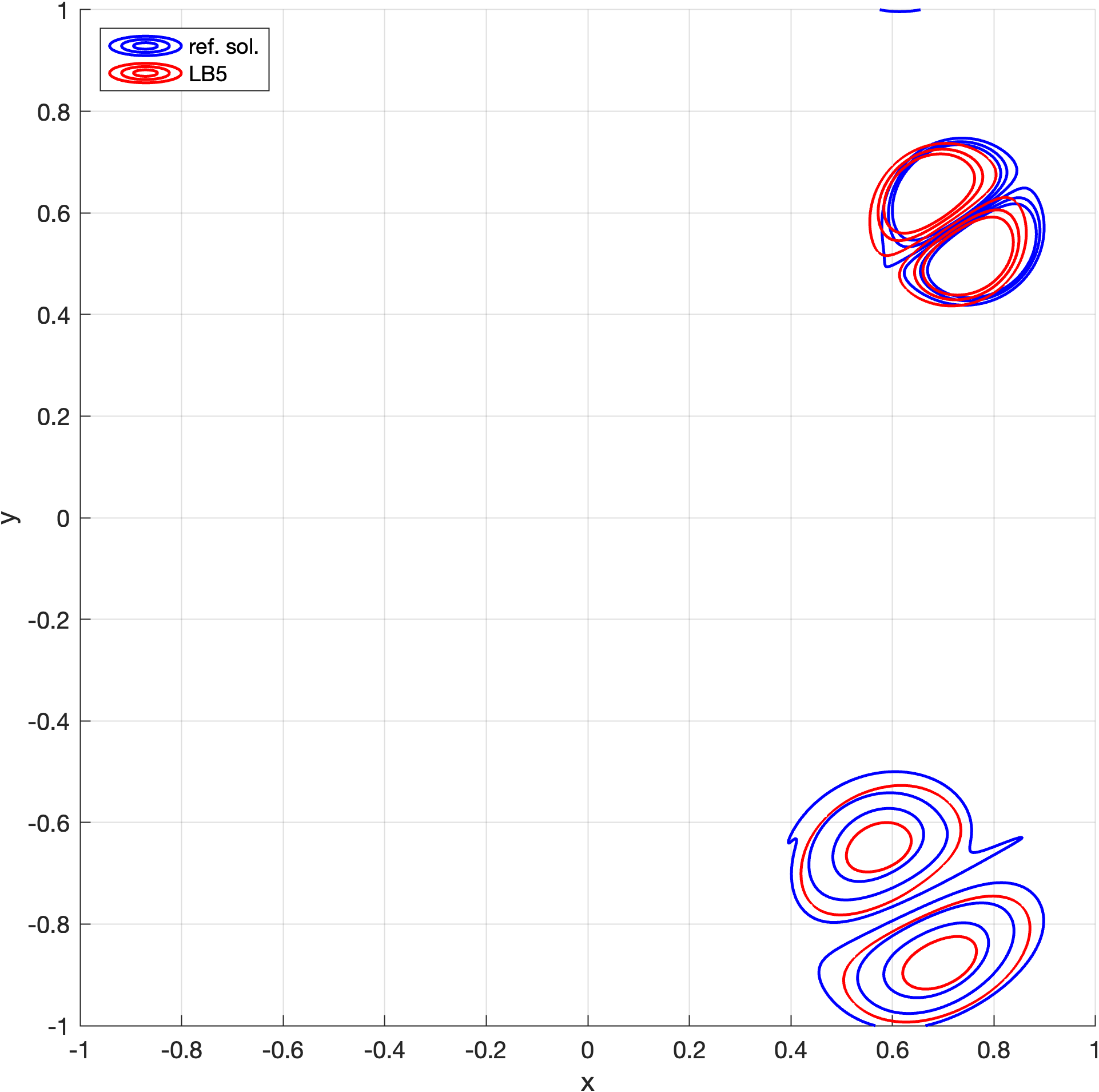}
        \\ \small (b) Spectral vs. LB5 contours
    \end{minipage}

    % Légende globale en anglais pour votre article
    \caption{Vorticity field results at $t = 1.8$: (a) computed vorticity field using the LB5 scheme, and (b) comparison of the vorticity contours between the pseudo-spectral reference solution (blue) and the LB5 scheme (red).}
    \label{fig:comp_vor_LB5}
\end{figure}

\noindent The LB5 scheme is identical to the LB4 scheme, incorporating the Dubois equilibrium,
but with the free relaxation rates specifically chosen to achieve an isotropic scheme, as detailed in the previous section (see \cite{ADGL14} for further details).
Consequently, all relaxation rates are fixed as follows: $s_x$ is determined to achieve a Reynolds number of $Re=2500$, $s_\varepsilon = s_x$, $s_q$ is selected such that $\sigma_q \sigma_x = \frac{1}{6}$, and $s_h = s_x$.
Figure~\ref{fig:comp_LB5} compares the kinetic energy, enstrophy, and palinstrophy computed by the LB5 scheme with those of the reference solution.
%The results demonstrate a good agreement between the LB1 simulation and the reference solution.
Figure~\ref{fig:comp_vor_LB5}(a) displays the vorticity field computed by the scheme (LB5) at time $t=1.8$.
Figure~\ref{fig:comp_vor_LB5}(b) compares the vorticity contours of the LB5 solution with those of the reference solution.
The figure clearly shows that the positions of the vortices obtained with the LB5 method are very close to those of the reference solution.
This indicates that the LB5 scheme exhibits the smallest anisotropy error compared to the other LB schemes investigated in this study.
Consequently, it is the model that best preserves the Galilean invariance property.

\section*{Conclusion}

\noindent In this work, a comprehensive numerical assessment of five lattice Boltzmann (LB) schemes was conducted to evaluate their accuracy regarding isotropy, numerical stability, and Galilean invariance. The highly convective oblique vortex dipole benchmark at $Re = 2500$ was utilized as a rigorous diagnostic tool, with numerical results systematically validated against a high-resolution pseudo-spectral reference solution.

The comparative analysis yielded several critical insights regarding the interplay between moment equilibrium formulations and relaxation rates. First, while both the standard Multiple-Relaxation-Time model (LB1) and the projected LBM formulation (LB3) successfully capture global integral quantities, they exhibit a distinct spatial phase lag in the positioning of vortical structures. This spatial error is directly attributable to a severe lack of isotropy and the presence of spurious terms that break Galilean invariance. Nevertheless, the projected formulation (LB3) demonstrated a clear operational advantage in terms of numerical stability and enabling robust computations at lower bulk viscosities where the standard MRT scheme undergoes catastrophic divergence.

Second, the investigation underscores that structural improvements to the equilibrium distribution function must be paired with careful parameter calibration to fully restore isotropy. Integrating the Dubois equilibrium for the heat flux moments (LB4) curtails the magnitude of the spurious terms but fails to completely rectify the directional bias. True multi-directional isotropy was achieved exclusively by the LB5 scheme, wherein the Dubois equilibrium is coupled with a precise tuning of the free relaxation rates ($\sigma_q \sigma_x = 1/6$). The LB5 model successfully minimized anisotropy errors, accurately preserving both the trajectories and the internal topology of the co-rotating vortex cores, thereby satisfying Galilean invariance.

Finally, the severe failure of the LB2 scheme—induced by truncating the non-linear contributions of the fourth-order moment equilibrium ($h^{eq}=0$)—provides important evidence regarding higher-order consistency. Although this truncation is known to have a negligible impact on standard, non-advective benchmarks such as Poiseuille or Taylor--Green flows, it triggered unphysical dissipation and the rapid destruction of the dipole in the present convective configuration. This finding demonstrates that the full non-linear equilibrium distribution function originally introduced in~\cite{QHL92} remains an absolute prerequisite for maintaining physical fidelity and preventing Galilean invariance violations beyond second-order hydrodynamic accuracy.

In summary, the oblique dipole benchmark has proven to be an exceptionally sensitive framework for isolating discrete directional errors. The results established herein suggest that advancing the predictive capability of lattice Boltzmann methods for complex turbulent or convective flows requires not only rigorous higher-order equilibrium definitions but also a systematic optimization of the relaxation parameters to strictly enforce isotropy and Galilean invariance across all resolved scales.

\bibliographystyle{plain} % Ou "alpha", "unsrt" selon le style désiré
\bibliography{bibliographie}
\end{document}